\documentclass[12pt]{amsart}
\usepackage{latexsym,fancyhdr,amssymb,color,amsmath,amsthm,graphicx,listings}
\usepackage[section]{placeins}
\definecolor{red1}{rgb}{1,0.9,0.9} \definecolor{blue1}{rgb}{0.9,0.9,1} \definecolor{green1}{rgb}{0.9,1,0.9} 
\definecolor{yellow1}{rgb}{1,1,0.9} \definecolor{yellow2}{rgb}{1,1,0.8}

\def\question#1{ \vspace{2mm} \begin{center} \fcolorbox{green1}{green1}{ \parbox{11.2cm}{{\bf Question:} #1}} \vspace{2mm} \end{center} }
\def\conjecture#1{ \vspace{2mm} \begin{center} \fcolorbox{green1}{green1}{ \parbox{11.2cm}{{\bf Conjecture:} #1}} \vspace{2mm} \end{center} }
\def\remark#1{ \vspace{2mm} \begin{center} \fcolorbox{yellow1}{yellow1}{ \parbox{11.2cm}{{\bf Remark:} #1}} \vspace{2mm} \end{center} }

\title{Goldbach for Gaussian, Hurwitz, Octavian and Eisenstein primes}
\author{Oliver Knill}
\date{Jun 19, 2016}
\address{Department of Mathematics \\ Harvard University \\ Cambridge, MA, 02138 }

\begin{document}
\maketitle

\begin{abstract}
We formulate Goldbach type questions for Gaussian, Hurwitz, Octavian
and Eisenstein primes. 
\end{abstract}

\section{Introduction}

The {\bf Goldbach conjecture} \cite{Guy} stating that every even integer $>2$ can be written as a
sum of two rational primes has a calculus reformulation in 
that the function $g(x)=f(x)^2$, with $f(x) = \sum_p x^p/p!$ summing over all primes has positive 
derivatives $g^{(2k)}$ for every $k>1$. Analogue questions can be asked in the other division 
algebras $\mathbb{C},\mathbb{H},\mathbb{O}$ as well as in number fields within $\mathbb{C}$. 
There is a Goldbach version by Takayoshi Mitsui for rings of integers \cite{Mitsui} and a
Goldbach conjecture by C.A. Holben and James Jordan for Gaussian primes
\cite{HolbenJordan}. Guided by the calculus reformulations, we 
look at Gaussian, Eisenstein, Quaternion and Octonion versions and make them plausible 
by relating them to conjectures by Edmund Landau, Viktor Bunyakovsky and Godfrey Hardy and John Littlewood. 
Even before the 1742 letter from Christian Goldbach to Leonard Euler, \cite{HelmutKoch,GoldbachBuch}. 
Ren\'e Descartes voiced a similar conjecture earlier \cite{dicksonI,Schechter,Vaughan,Pintz}:
every even number is the sum of 1,2 or 3 primes.
Edward Waring conjectured in 1770 that every odd number is either a prime or the 
sum of three primes. Many have done experiments. Even the Georg 
Cantor (1845-1918) checked in 1894 Goldbach personally up to 1000 \cite{dicksonI} p. 422.
The ternary conjecture is now proven  \cite{Helfgott2014}. It required to search up to
$n \leq 8.875 10^{30}$. The binary has been verified to $4 \cdot 10^{18}$ \cite{OHP}. 
Landmarks in theory were Hardy-Littlewood \cite{HardyLittlewood1923} with the circle method, 
with genesis paper with Srinivasa Ramanujan in 1918 \cite{HardyRamanujan,Vaughanbook},
the Lev Schnirelemans theorem \cite{Hua1982}
using density and Ivan Vinogradov's theorem \cite{Vinogradov} using trigonometric sums.

\section{Gaussian primes}

A Gaussian integer $z=a+ib$ is a {\bf Gaussian prime} if it has
no integer factor $w$ with integer norm $N(z)=a^2+b^2=|z|^2$ satisfying
$1<|w|^2<|z|^2$. A Gaussian integer is prime if and only if $p=N(z)=a^2+b^2$ is a 
rational prime or if $ab=0$ and $\sqrt{p}$ is a rational prime of the form $4k-1$. 
This structure, which relies on Fermat's two square theorem (see \cite{Zagier90} for a topological proof),
has been known since Carl Friedrich Gauss \cite{Gauss1831}:
there are three type of primes: the case $\pm 1 \pm i$ for $p=2$, then primes of the 
form $\pm p$ or $\pm i p$ for rational primes $p=4k+3$ or then groups $\pm a \pm i b, \pm b + i \pm a$ of 
eight primes for rational primes of the form $p=4k+1$. \\

Lets call $Q=\{ a+ib \; | \; a>0, b>0 \}$ the {\bf open quadrant} in the complex integer plane. 
A Gaussian integer $z$ is called {\bf even} if $N(z)$ is even. Evenness is equivalent to
$a+b$ being even or then of being divisible by the prime $1+i$. 

\begin{figure}[!htpb]
\scalebox{0.4}{\includegraphics{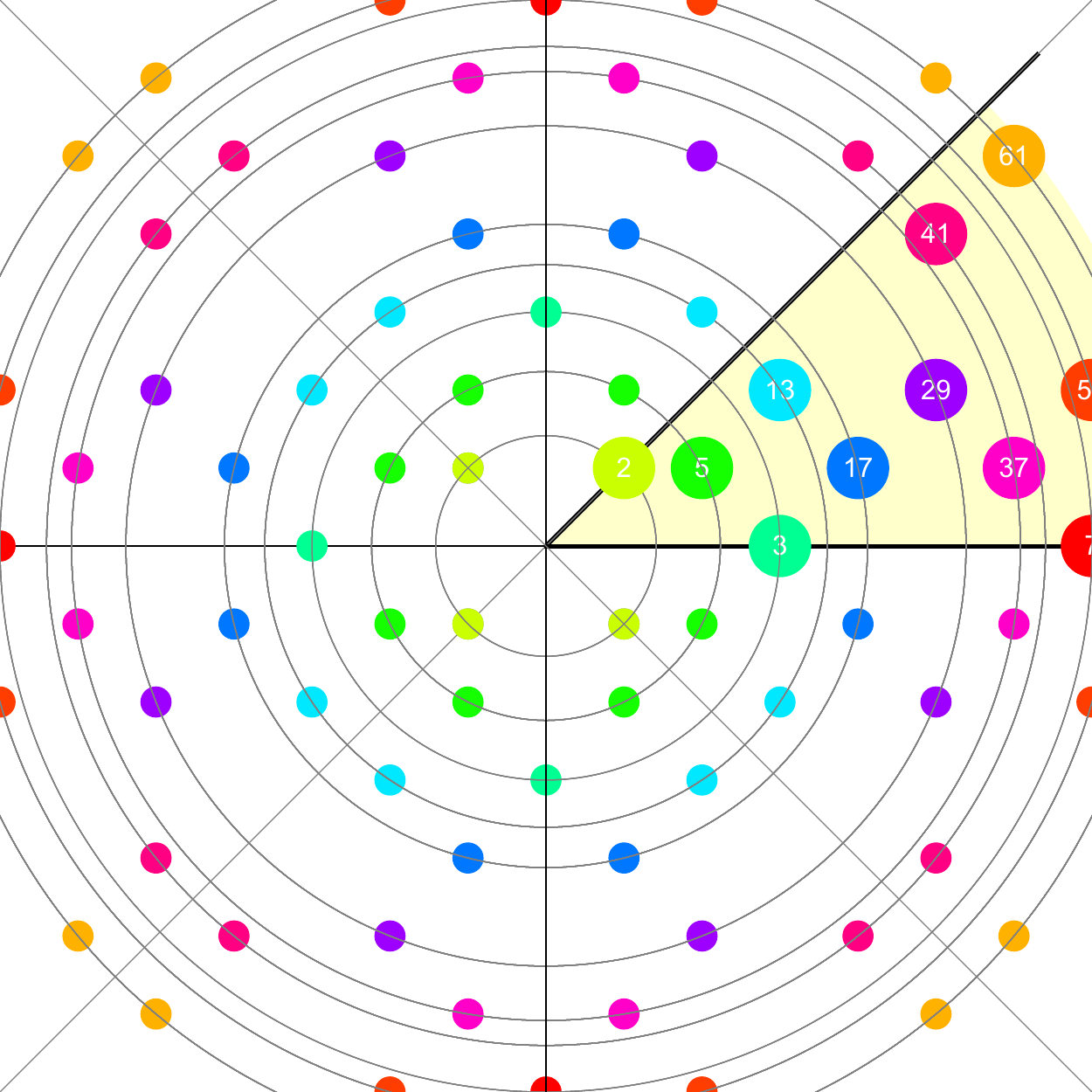}}
\scalebox{0.4}{\includegraphics{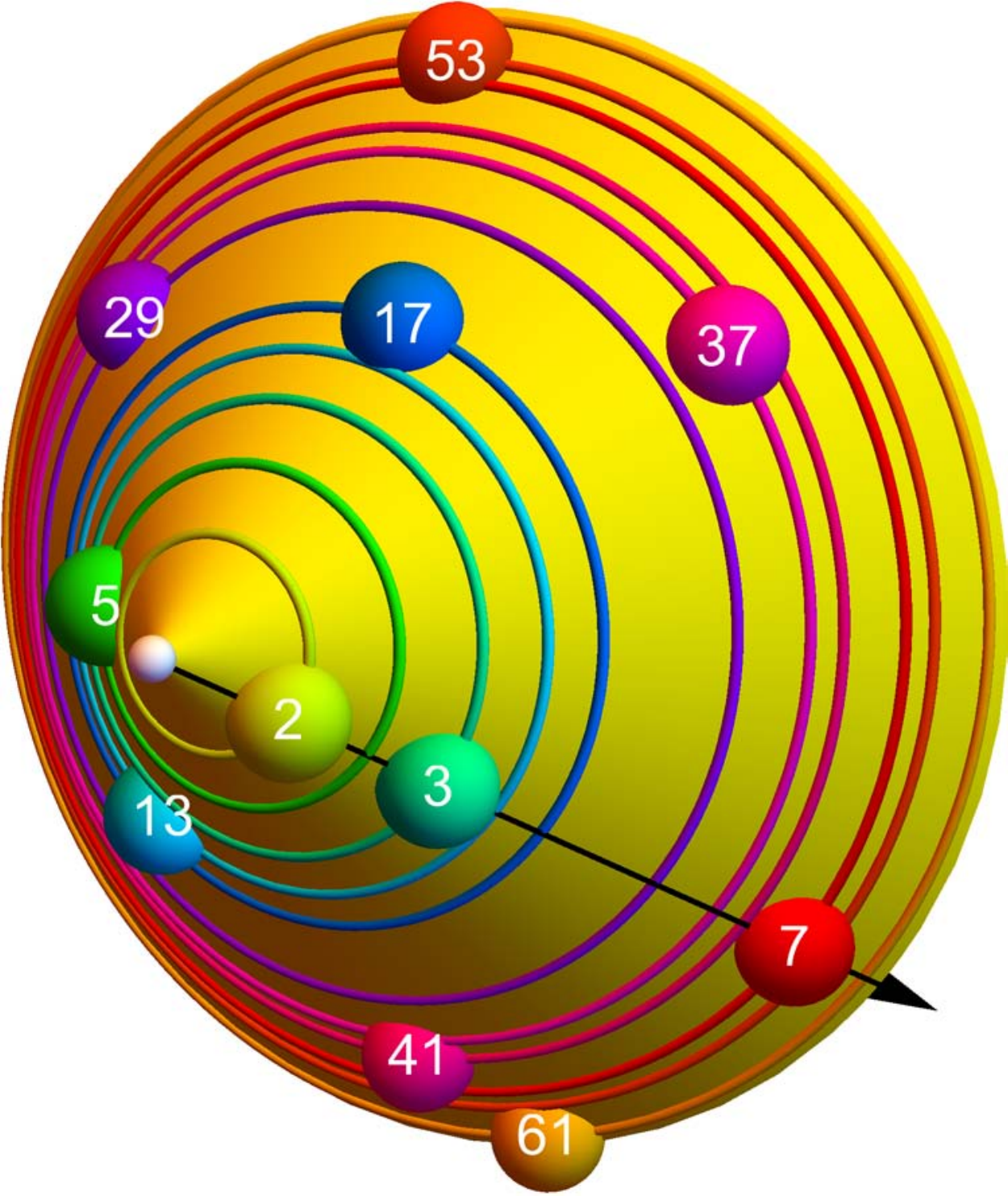}}
\caption{
Gaussian primes cover the rational primes in a natural way. Only the order is scrambled
although. There is a dihedral $D_4$ symmetry of the primes generated by conjugation and
multiplication by units. On the cone $\mathbb{C}/D_4$, the angle distribution 
appears pretty random away from the identification line.
}
\end{figure}

\conjecture{
Every even Gaussian integer $a+ib$ with $a>1,b>1$ in $Q$ is the sum of two Gaussian primes in $Q$.
}

This can again be reformulated using Taylor or Fourier series. In the Taylor case, the conjecture is equivalent 
with 
$$  f(x,y) = \sum_{p=a+ib \in Q} \frac{x^a}{a!} \frac{y^b}{b!} $$ 
having the property that $g=f^2$ has all even derivatives 
$g^{(k,l)}$ positive if $k+l$ is even and $k>1,l>1$. The Fourier case appear with 
$x=\exp(i \alpha), y = \exp(i \beta)$. Such algebraic reformulations make the statement natural. 
The conjecture looks toughest at the boundary of $Q$, where less possibilities for summation appear. The extremest
case is $z=2n+2i$, where $z=(a+i) + (b+i)$ forces $2n=a+b$ with $a^2+1,b^2+1$ both being prime. We see:

\remark{
Gaussian Goldbach implies the Landau conjecture on the infinitude of primes of the form $n^2+1$.
}

Since Landau appears currently out of reach, proving the Gaussian Goldbach will not be easy. 
There is still the possibility although of a counter example.
It looks unlikely, given the amazing statistical regularity predicted by Hardy and Littlewood. 
But a surprise is always possible. Here is the ternary version which like in the real case 
does not require an evenness condition:

\conjecture{Every Gaussian integer $a+ib$ in $Q$ satisfying $a>2,b>2$ 
is the sum of three Gaussian primes in $Q$.}

Lets compare with what has been asked before: \\
Holben and Jordan \cite{HolbenJordan} formulate as `conjecture D" the statement that every even Gaussian integer is the 
sum of two Gaussian primes, and then ``conjecture E": if $n$ is a Gaussian integer 
with $N(n)>2$, then it can be written as a
sum $n=p+q$ of two primes for which the angles between $n$ and $p$ as well as $n$ and $q$ are
both $\leq \pi/4$. Their conjecture F claims that for $N(n)>10$, one can write $n$ as a sum of two primes
$p,q$ for which the angles between $n$ and $p$ and $n$ and $q$ are both $\leq \pi/6$. \\
Mitsui \cite{Mitsui} formulates a conjecture for number fields which of course include
Gaussian and Eisenstein cases. In the Gaussian case, it states that every even Gaussian integer is the 
sum of two Gaussian primes. This is conjecture D in \cite{HolbenJordan}.
Even when using the full set of primes, the evenness condition is necessary. 
The smallest Gaussian integer which is not the sum of two Gaussian primes is $4+13i$. The smallest real one is $29$. 
The Holben-Jordan conjecture implies the Mitsui statement.  In the Eisenstein case, we see
that every Eisenstein integer is the sum of two Eisenstein primes without evenness condition. 
The question makes sense also in the $\mathbb{Z}$: is every even integer the sum of two {\bf signed primes},
where the set of {\bf signed primes} is $\{ \dots, -7,-5,-3,-2,2,3,5,7, \dots \}$. The smallest number
which is not the sum of two signed primes is $23$. By Goldbach, all even numbers should be the sum 
of two signed primes. The sequence of numbers 
\cite{A110673}. 
which are not representable as the sum of two signed primes 
contains the number 29 as the later is also not a sum of two Gaussian primes. 
We will predict however that every rational integer to be the sum of two Eisenstein primes simply because every 
Eisenstein integer seems to be the sum of two Eisenstein primes. 

\begin{figure}[!htpb]
\scalebox{0.9}{\includegraphics{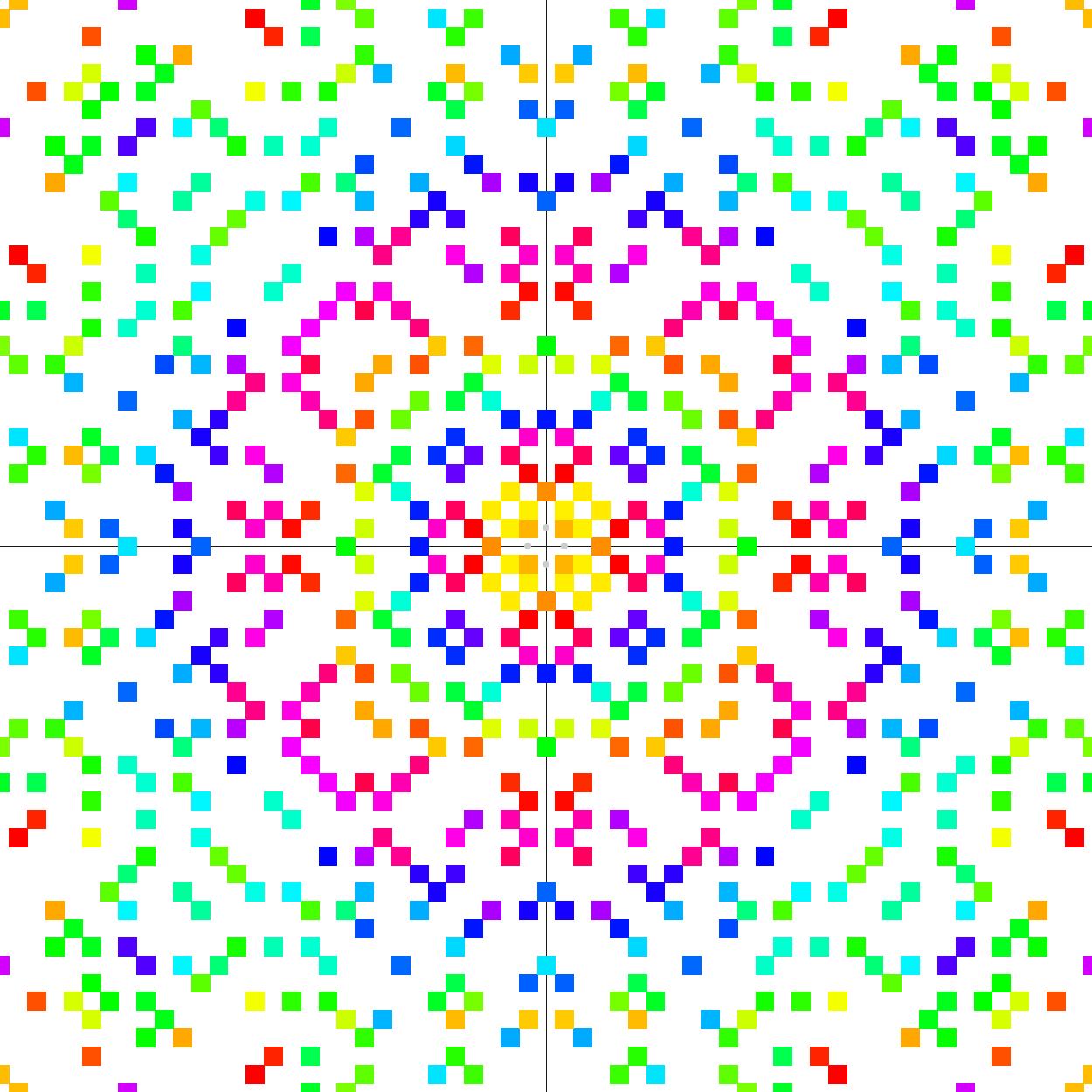}}
\caption{
Gaussian primes.
}
\end{figure}

\section{A Hardy-Littlewood constant}

The statistics of Gaussian primes on various rows has been of interest for almost 100 years, sometimes 
without addressing the Gaussian primes although. It is related to the fascinating story of a constant: 
{\bf Hardy-Littlewood} conjectured that the frequency ratio of primes on 
${\rm im}(z)=1$ and ${\rm im}(z)=0$ is 
$$ C = C_1= \prod_{p \in P_1} [1-\frac{1}{p-1}] \prod_{p \in P_3} [1+\frac{1}{p-1}] = 1.37279 ... \; , $$
where $P_k$ is the set of rational primes congruent to $k$ modulo $4$. More generally, 
in their conjecture $H$, they give explicit constants for all the other rations 
between rows $k+ia$ and $k+ib$ as $C_a/C_b$ using the product 
$C_a = \prod_{p} (1-(-a^2|p)/(p-1))$ over all odd primes $p$, where
$(q|p)$ is the {\bf Jacobi symbol}, assuming the empty product is $C_0=1$. \\

In 1962, \cite{BatemanHorn} combine several conjectures of Hardy and Littlewood about
density relations of sets. What is the ratio of the number of primes of 
the form $f(x)$ in $[0,\dots,n]$ and the number of primes the form $g(x)$ in $[0,\dots, n]$?
We can write this in the limit $n \to \infty$ as $C_f/C_g$, where
$$   C_f = \prod_p \frac{(1-\omega_f(p)/p)}{(1-1/p)} \; . $$
This intuition works if $f,g$ are irreducible polynomials of the same degree
with positive leading coefficients and $\omega_f(p) \neq p$ for all $p$
where $\omega_f(p)$ is the number of solutions $f(x)=0$ modulo $p$. 
One can include $f(x)=x^n$, where $\omega_f(p)=1$ and so set $C_f=1$.  \\

The constants $C_k$ are then shorts cuts for $C_{n^2+k}$. But lets just focus on the constant $C$ 
for $f(x)=x^2+1$ which is a nice prototype as in that case $\omega_f(p)=1$ for $p \in P_3$ and 
$\omega_f(p)=2$ for $p \in P_1$ by quadratic reciprocity so that 
$(1-\omega_f(p)/p)/(1-1/p)) = 1-\frac{1}{p-1}$ for $p \in P_1$ 
and $(1+\frac{1}{p-1})$ for $p \in P_3$. \\

The intuition behind the constants $C_f$ is of probabilistic nature. It is
a bit hard to describe but once you see it, 
the ``formulas open up" and become crystal clear. Lets try: 
we are interested in solutions $a^2=-1$ modulo $p$ because then, 
$p$ is a factor of $a^2+1$ preventing the later of being prime. 
The key assumption is that the solution sets of equations like $a^2=-1$ modulo $p$ or then modulo $q$ 
is pretty much independent events if $p,q$ are different odd primes. This means that the chance being a solution
to both becomes the product of the probabilities for each. To get the probability, start with the full set of
integers on some large interval, then look at the probability changes, when primes are successively added to the list.
Now, every time a new prime is added to the list, the size of the space changes by 
$(p-1)/p=(1-1/p)$ because we can only take numbers which are not multiples of $p$. The product of these size 
changes is $1/\zeta(1)=0$ and reflects the infinitude of primes. But if we look at the {\bf ratio} 
of solution sets for two polynomials $f,g$, we don't have to worry about this renormalization: 
it happens for both $f$ and $g$ in the same way. 
In particular, when looking at solutions of the form $x^2+1$, then whenever $-1$ is 
a quadratic residue, the probability decreases by $(1-1/(p-1)) (1-1/p) =(p-2)/p$ but if $-1$ is 
not a quadratic residue, then the probability stays the same, as the prime has no chance of dividing $a^2+1$.
Including again the volume change gives $(p-2)/p /(1-1/p) = 1-1/(p-1)$ in the residue case and 
$1/(1-1/p) = p/(p-1) = 1+1/(p-1)$ in the non-residue case. This explains the formula
for $C$. The more general case is combined skillfully in \cite{BatemanHorn}. Using this frame work,
many of the formulas of \cite{HardyLittlewood1923} make sense, also density formulas for the estimated number of
prime twins, where one takes $f(x)=x(x+2)$ as this is also a special case of the Paul Bateman and Roger 
Horn setup. This general statement generalizes many conjectures and is called the 
{\bf Bateman-Horn} conjecture. 
For relations with Dirichlet series see \cite{KConrad2003} for generalizations
to number fields see \cite{GrossSmith}. For computations of various constants related to patterns 
of Gaussian integers in particular, see \cite{GaussianZoo}.  \\

Mathematicians were interested in these constants also for
cryptological reasons. A major reason is because there is a {\bf holy grail} 
for factoring large integers $n=pq$ known already to Pierre de Fermat:
find solutions to quadratic equations like $x^2=a$ modulo $pq$ for small $a$. Since $x^2-y^2$ is then zero modulo $pq$
the greatest common divisor of $x-y$ and $n$ is a factor. Essentially all advanced factorizaton methods 
like Morrison-Brillard or quadratic sieve are based on this \cite{Riesel}. Also in this setup, some heuristic
randomness assumptions for pseudo random sequences are the key to estimate the time an algorithm needs to factor the
integer. There were later computations by Daniel Shanks \cite{Shanks60} and Marvin Wunderlich (1937-2013) 
\cite{Wunderlich} who was also a cryptologist.
In 1973, Wunderlich checked the frequency ratio of primes up to 14'000'000 and compared it with 
the Hardy-Littlewood constant $C$. We profited from Moore's law on hardware advancement and computed 
up to $134143000000$  
up to which $2728969165$ rational primes of the form $4k+3$ exist and $3746378328$ primes of the form $n^2+1$
which gives a fraction of $1.372818123..$. Fluctuations still happen.

\begin{figure}[!htpb]
\scalebox{0.3}{\includegraphics{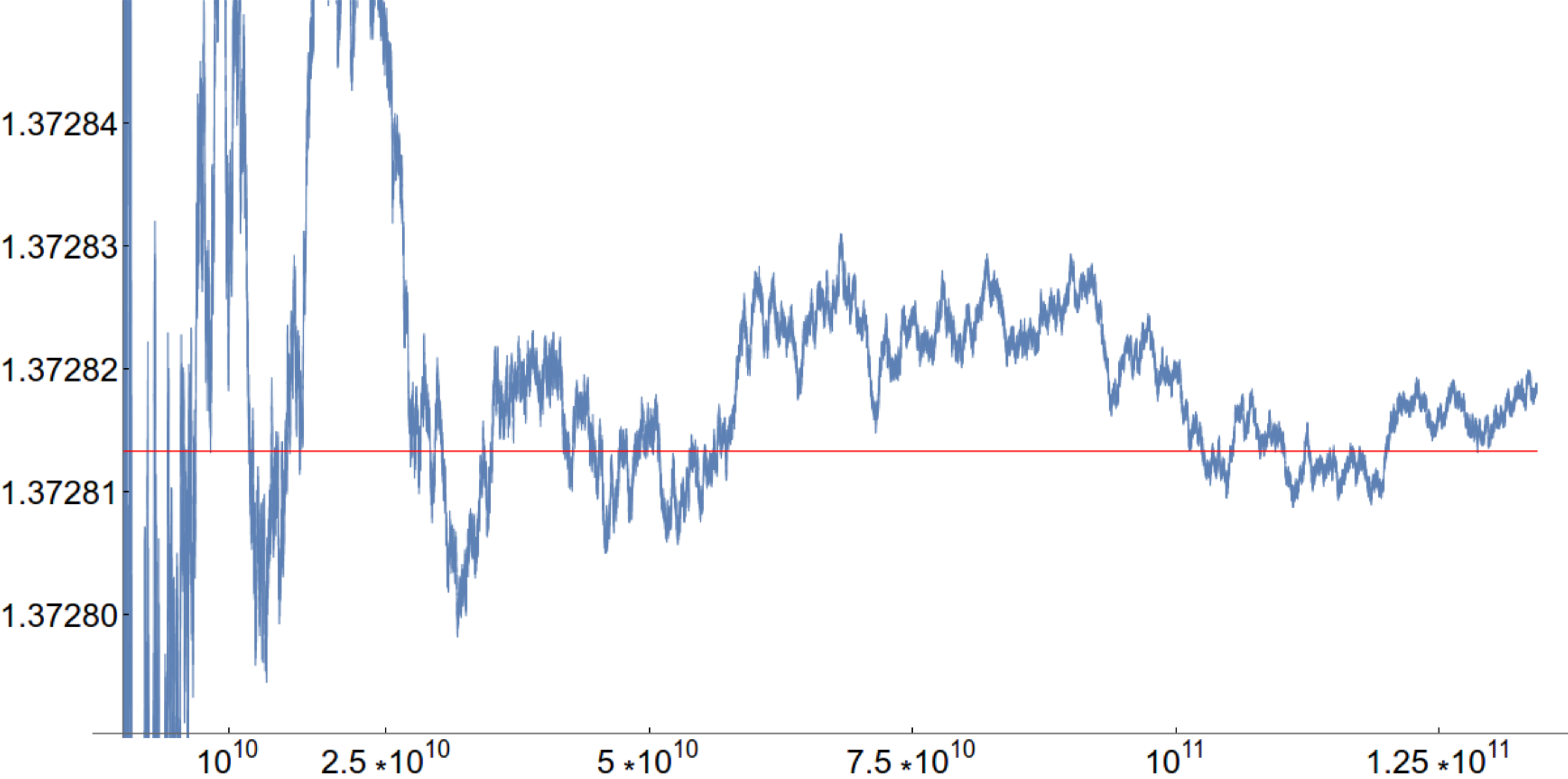}}
\caption{
The convergence to the Hardy-Littlewood constant $C$ giving the 
ratio of Gaussian prime density on ${\rm Im}(z)=1$ and ${\rm Im}(z)=0$. 
Shanks \cite{Shanks60} gave $C=1.37281346$.
The constant $C$ is almost prophetic as by an open Landau's problem (which currently
appears theoretically beyond reach), one does not even know whether $C$ is positive!
}
\end{figure}

By the way, the grandmasters of number theory, Godfrey Hardy and John Littlewood were not without assistance: 
we read on page 62 of \cite{HardyLittlewood1923} that "some of their conjectures have been tested
numerically by Mrs. Streatfield, Dr. A.E. Western and Mr. O. Western". Much has been written about the
influence of computers in mathematical research \cite{Williams82}. The story of the Hardy-Littlewood 
constant $C$ illustrates that already early in the 20th century, when humans were still doing the computations by hand,
the experimental part has been important. We were especially intrigued by the constant $C$ 
because Streatfield, Western, Western already got a 5 digit agreement. This was later in the century confirmed by 
Shanks and Wunderlich. As we have much faster computers now, how accurate can we get now? Surprisingly, 
even for $n=2^{36}$, the fluctuations remain of the order $10^{-5}$. Surprised, we started to 
believe initially that the Hardy-Littlewood conjecture could be too strong and some small fluctuations could remain.
There is no reason for alarm however. Theory predicts a slow convergence: if the Riemann hypothesis holds, an error of the
order $(\log(x))^2/\sqrt{x}$ is expected. While this converges to zero, it happens very slowly in the range 
we do currently compute: for $x=2^{36}$, this is still $2 \cdot 10^{-3}$ only. So, the constant is already in 
better agreement than expected. Just to get an idea about the computing culture:
Daniel Shanks (1917-1996) \cite{Shanks60} reports in 1959 that A.E. Western rewrote the constant using
zeta and beta function values as 
$$ C = \frac{6}{4} \frac{\zeta(6)}{\beta(2) \zeta(3)} 
       \prod_{p \in P_1} (1+\frac{2}{p^3-1}) (1-\frac{2}{p (p-1)^2}) $$ 
which gives 4 decimal places already when summing over three primes $p=5,13,17$: the value is then
$\frac{478543065 \pi ^6}{304368582656 \beta(2) \zeta(3)}=1.37283...$, 
where $\beta(2)=0.915966$ is the Catalan number and $\zeta(3)=1.20206$ is a zeta value. 
A bit earlier, in 1922, even before the Hardy-Littlewood article appeared,
A.E. Western \cite{Western22} computed the constant $C$ using further sophisticated identities
involving various zeta values so that one can use two primes $5,13$ only to get $C$
to 4 decimal places! Western was a giant in computation who also published mathematical tables like
\cite{WesternMiller}. Anyway, this episode illustrates the culture in which the almost prophetic
paper \cite{HardyLittlewood1923} was written in.  \\

Its maybe important to remember that the mathematicians one hundred years ago had no access to computers;
and this situation remained essentially until the mid century. Zuse Z3 was built 1941, Colossus and Mark I 
appeared only in 1944 but all of them were slow: 
Mark I for example needed 6 seconds to multiply two numbers \cite{NahinNumberCrunching}. 
Shanks in 1960 used an IBM 704 with a 32K high-speed memory. Today we have access to machines
which give each user 500 GBybes of RAM and hundreds of CPUs.
The program of Shanks needed on the Vaccuum tube computer 
10 minutes to factor all $n^2+1$ from $n=1$ to $n=180000$. Today, a 
tiny laptop weighting 2 pounds with 8 GBytes of RAM reports it done in 
less than 4 seconds. 

\section{More questions}

Much about the structure of Gaussian integers is still a mystery. One does not know for example
whether there are some rows without primes. This is the {\bf frogger problem}: the complex plane is the freeway,
the rows are the car lanes, the non-primes are the cars. Can the frog hop to infinity on primes, which are the car gaps?
The {\bf Gaussian moat problem} has attracted a lot of attention: does a bounded step size
suffices to hop to infinity on primes?  
\cite{JordanRabung,GethnerStark,GethnerWagonWick,Wagon,Loh,Vardi}.
The {\bf twin prime conjecture} for Gaussian primes
asks for the existence of infinitely many {\bf Gaussian primes twins}, pairs of primes for
which the Euclidean distance is $\sqrt{2}$ \cite{HolbenJordan}.
While one does not know whether infinitely many Gaussian prime twins exist, one
can estimate that there are asymptotically $C r/\log^2(r)$ of them in a ball of
radius $r$ \cite{GaussianZoo} which is a Hardy-Littlewood type estimate and now part of the Bateman-Horn 
conjecture. Virtually any question known for rational integers can be ported to Gaussian primes. 
This includes Diophantine problems or Waring type problems. Here is an other question: what is the regularity of the sequence
defined by the unique angle $\theta(p_k) \in (0,\pi/4)$ if $p_k$ is the
$k$'th Gaussian prime in that sector? The sequence behaves very much like a random number generator.
Is the topological entropy positive for the shift on the compact hull generated by the sequence? \\

While the numerical evidence for Hardy-Littlewood is strong, one 
has to get reminded that even the existence of infinitely many primes on the Gaussian line $\{ {\rm Im}(z)=1\}$ 
is open. It is one of the four problems presented by {\bf Edmund Landau} at the 1912 International congress of mathematicians.
The statement beats even Goldbach in simplicity: are there infinitely many primes $p$ for which $p-1$ is a square?
Hardy and Littlewood have also formulas for the density of Gaussian primes of the form $\{ a + i b_0 \}$ for fixed $b_0$. 
We looked experimentally at correlations between primes in different rows or at matrices and constructed
graphs defined by Gaussian primes: take as the vertex set the positive integers and connect $a,b$ if
$a+ib$ is a Gaussian prime. We also applied the game of life to Gaussian prime configurations and believe
that {\bf Gaussian life} exists arbitrary far away from the origin. It could be ``prime twins" blinking
from far, far away. An other question is if we look at Gaussian or Eisenstein primes on the orbifold
$\mathbb{C}/D_4$ or $\mathbb{C}/D_6$ factoring out the symmetry and let the ones with integer norm (not the ones with 
integer radius) move with uniform speed on those
cones (like planets circling the sun). Is there any positive time for which 4 primes $a+ib$ with prime $a^^2+b^2$
are on the same line?  It is likely that such {\bf exceptional prime alignments} do not exist for Gaussian 
primes and similarly for Eisenstein primes away from the symmetry axes. 
Also, when factoring out the symmetry group given by the units and conjugation in the Hurwitz or Octonion case,
one gets prime constellations on compact manifolds. In the Hurwitz case, there is a compact $3$ manifold which
is diffeomorphic to $M(p)=(\mathbb{H} \cap S_{p})/G$, where $S_{p}$ is the sphere in $\mathbb{R}^4$ 
of radius $p$ and $G$ is the group generated by the action of multiplying with a units and conjugation. 
$M(p)$ plays the role of the circles through the primes centered at the origin on the cone in the Gaussian case. 
Hurwitz showed that there are exactly $p+1$ primes on $M(p)$. We expect the distribution of this 
{\bf prime cloud} on on $M$ to converge weakly to an absolutely continuous volume 
measure on $M$, when the rational prime $p$ goes to infinity. \\ 

\begin{figure}[!htpb]
\scalebox{0.4}{\includegraphics{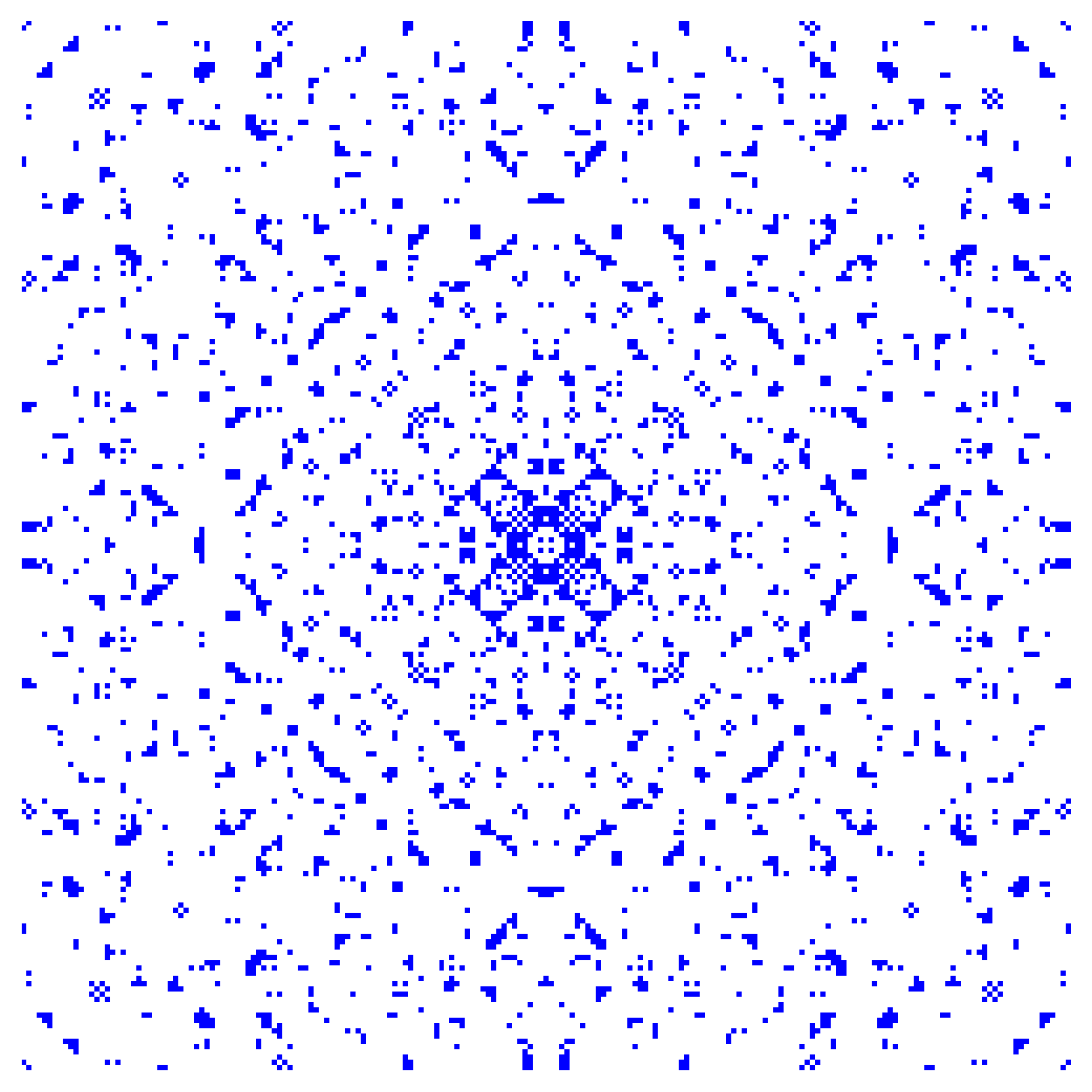}}
\scalebox{0.4}{\includegraphics{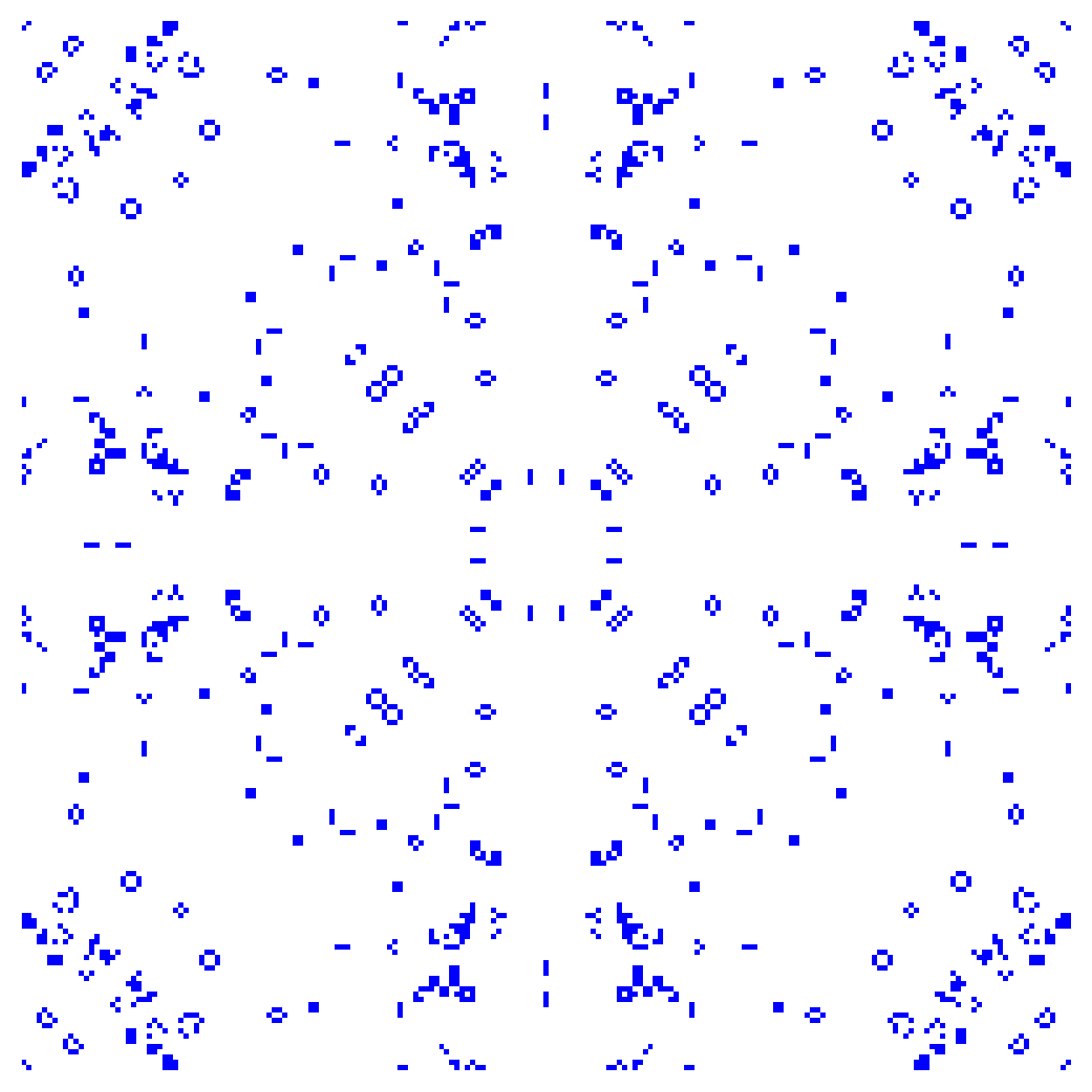}}
\caption{
Gaussian primes are a configuration in $A^{Z^2}$, where $A=\{0,1\}$ is the two letter
alphabet. We can therefore apply cellular automata on it. 
The left picture shows it after applying the Conway game of life rule
once, the second after three times. In light of the prime twin conjecture for
Gaussian primes, there should be life arbitrary far away from the origin. 
``Life" in a region is a configuration which ``moves" when applying the time map.
\label{life}
}
\end{figure}

But there are not only open problems about Gaussian primes. Quite many results are known.
There is an analogue of Dirichlet's theorem on arithmetic progressions:
for an arbitrary finite set in $\mathbb{Z}[i]$, there exist infinitely many
$a \in \mathbb{Z}[i]$ and $r \in \mathbb{R}^+$ such that $a+r \sum_{f \in F} v_f$
is a Gaussian prime \cite{Tao2006}. An other example is that the
density of the {\bf prime quotients} $p/q$ in $R^+$ generalizes to the statement
that the Gaussian prime quotients $p/q$ are dense in the complex plane \cite{Garcia}.
Patterns are explored in \cite{GethnerWagonWick,Wagon,JordanRabung76,GaussianZoo}. 
The Goldbach conjecture is not the only statement which involves the
additive structure and primes (which inherently rely on the multiplicative
structure of the ring): any additive function $f(zw)=f(z)+f(w)$ which satisfies
$f(p+1)=0$ for all Gaussian primes is $0$ \cite{MehtaViswandham}.
Gaussian primes and friends are an Eldorado for new questions. 
We got dragged into this area while doing an exercise in section 5 of 
\cite{MazurStein}. The addictive topic of primes totally ruined our original summer plans, just 
like it did for the boy in the Doxiadis novel \cite{petros}. 

\begin{figure}[!htpb]
\scalebox{0.1}{\includegraphics{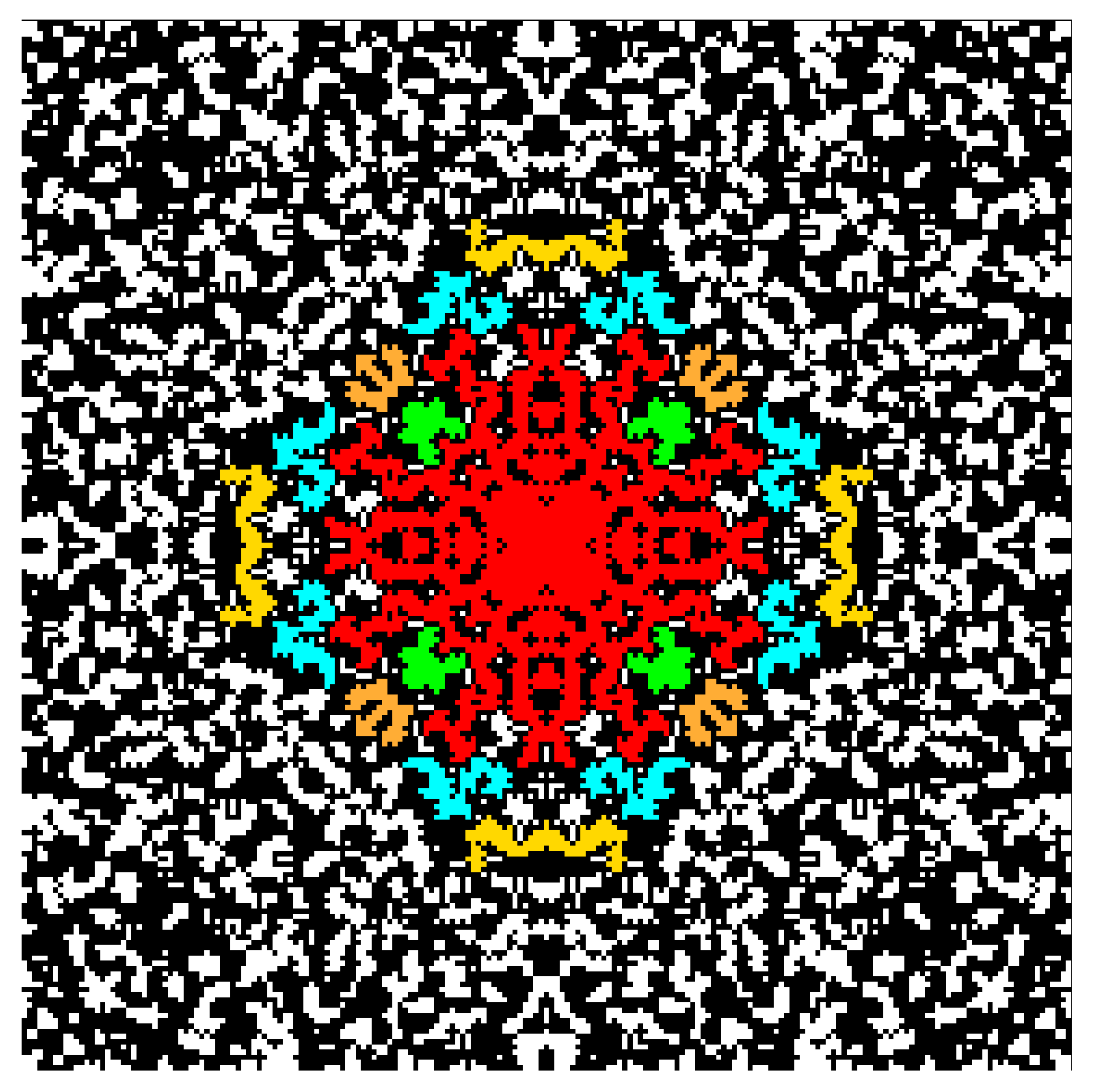}}
\scalebox{0.1}{\includegraphics{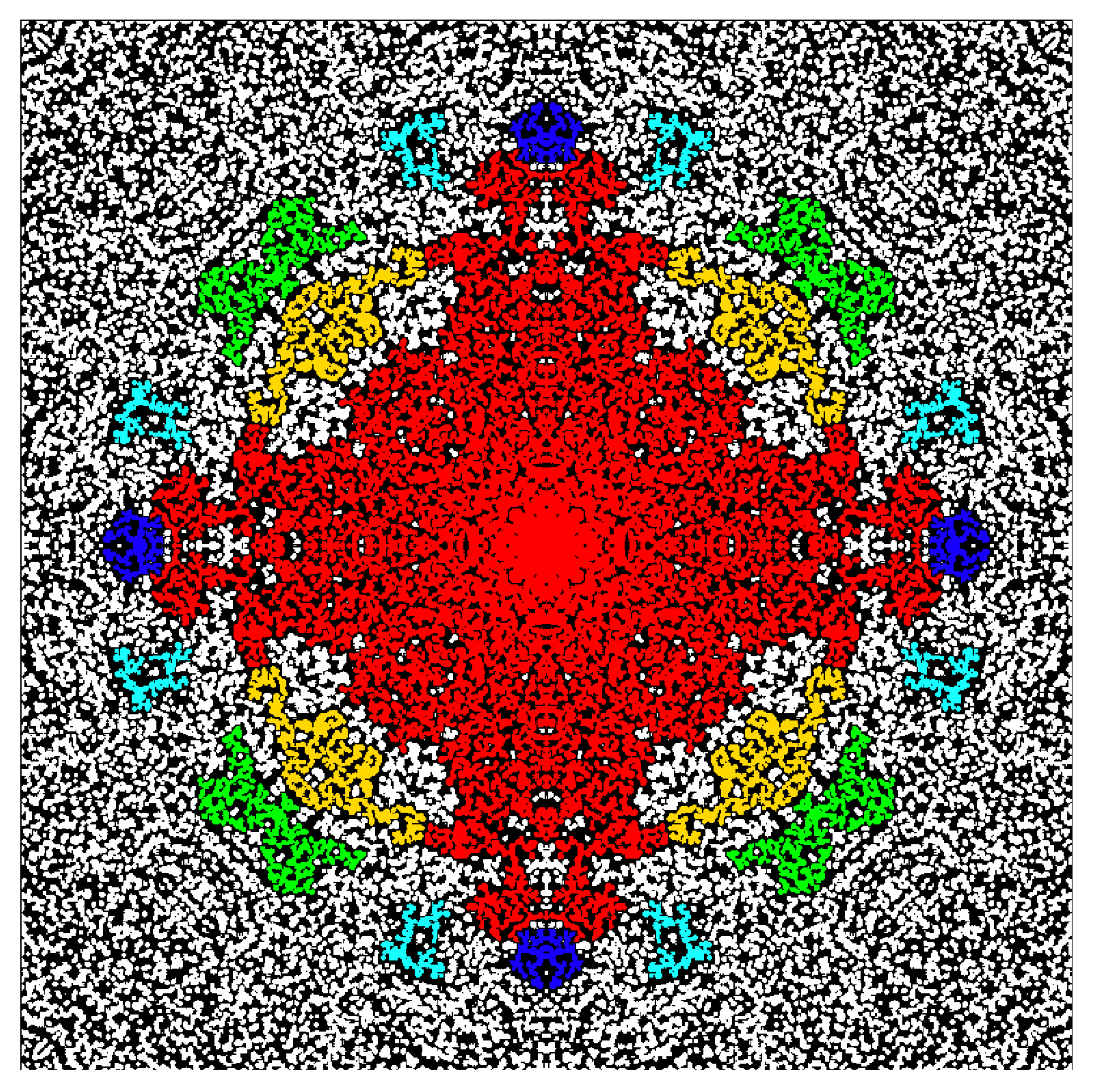}}
\caption{
To illustrate the moat problem we applied a cellular automaton
map and looked at the central connected component. The first
picture is a screen shot after applying the map once, the second after applying
the map two times. In principle, one can compute moats like that.
But it is not very efficient. 
\label{moat}
}
\end{figure}

\section{Eisenstein primes}

Given the cube root $w=(1+\sqrt{-3})/2$ of $-1$, an {\bf Eisenstein integer} is a complex
number $a+b w$ for which $a,b$ are rational integers. The norm $p=N(z)$ of $z=a+bw$ is $a^2+b^2+ab$. 
One usually writes the Eisenstein integers with the cube root of $1$. For Goldbach statements however, 
it is more convenient to work with $w$ which is in the first quadrant. 
This ring of integers has been investigated first by Gotthold Eisenstein (1823-1852) \cite{Eisenstein}.
The basic structure of {\bf Eisenstein primes} is well known \cite{StillwellNumber}. 
Either $p$ is prime and congruent $0,1$ modulo $3$ or then 
$\sqrt{p}$ is prime and congruent to $2$ modulo $3$. The later class are the primes on the 
six symmetry axes. \\

Define $Q$ as the set of Eisenstein integers $a+b w$ with $a>0,b>0$. These are the integers 
in the fundamental region, the first sextant of the complex plane. After doing some experiments, 
we experienced two surprises: first of all, the Eisenstein primes appear so dense that no 
{\bf evenness} condition is necessary for writing an integer as a sum of two primes. 
The second surprise was that the extreme boundary case appears less constraining than the 
nearest next boundary case. There are two isolated ``ghosts" Eisenstein primes.

\begin{figure}[!htpb]
\scalebox{0.4}{\includegraphics{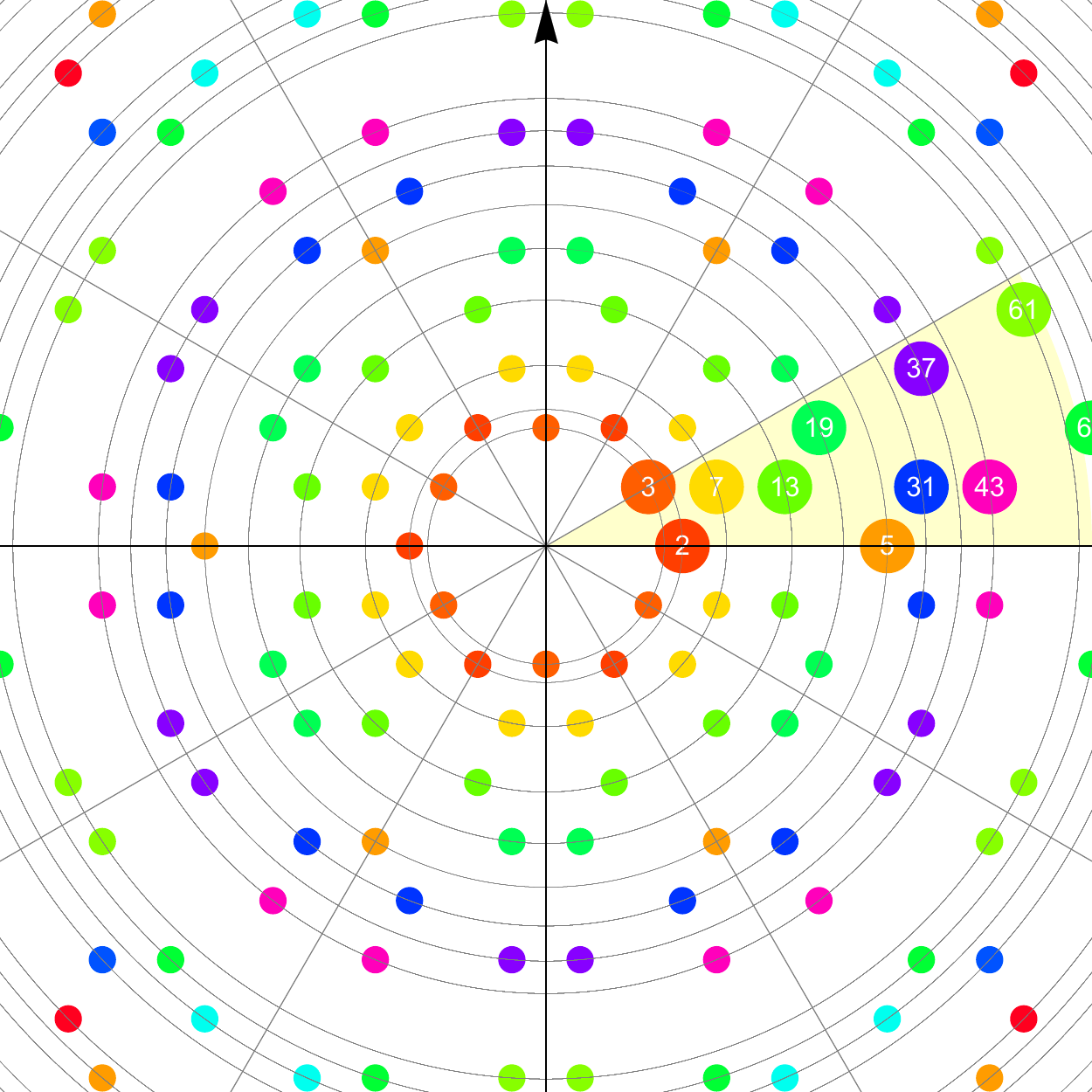}}
\scalebox{0.4}{\includegraphics{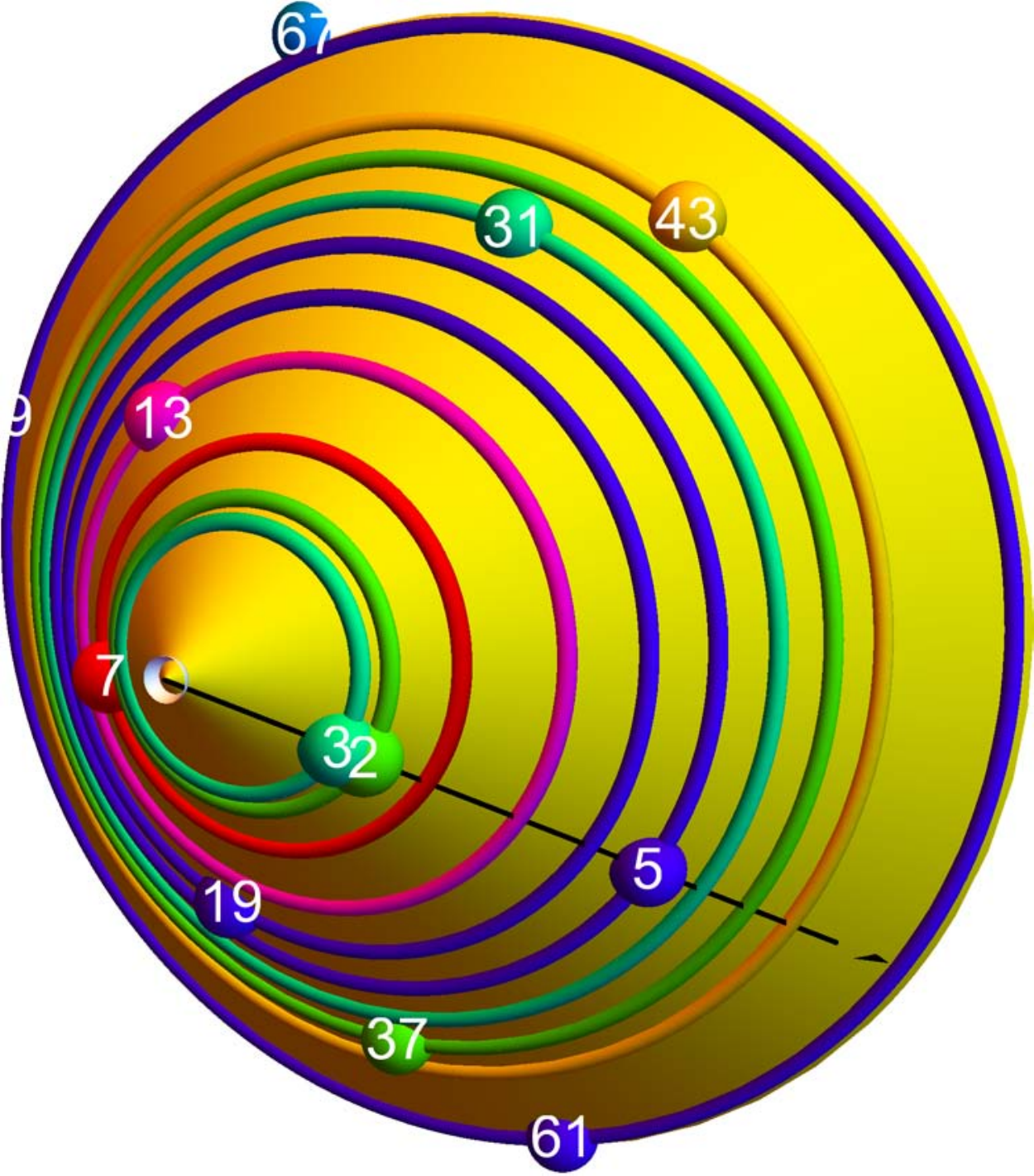}}
\caption{
Eisenstein primes also cover the rational primes. The order
is different as in the Gaussian case. The dihedral $D_6$ symmetry 
of the primes hides that on the fundamental domain
$\mathbb{C}/D_6$, the angle distribution is pretty random, apart from
the prime $3$ and the $3k+2$ primes confined to the real axes.
}
\end{figure}

\conjecture{
Every Eisenstein integer $a+bw$ with $a>3,b>3$ is the sum of two Eisenstein primes in $Q$. 
}

This would especially imply that infinitely many primes exist which are of one of the three 
forms $n^2+n+1,n^2+2n+4,n^2+3n+9$. Since in the second case, $n$ has to be odd so that $a=n+m$ is
even, we need for odd $a$ a decomposition in the form $(n+w) + (m+3w)$ which means:

\remark{
Eisenstein Goldbach implies the Bunyakowsky conjecture on the infinitude of primes of the form $n^2+n+1$. 
}

We originally got the impression that the condition $a>1,b>1$ works and not only $a>3, b>3$, as stated. 
Here is a tale of caution: as we looked at the toughest case $z=2+nw$, where the problem is to write an 
integer $n$ as a sum $n=a+b$ where $a^2+a+1,b^2+b+1$ are both prime, then this problem appears always
to have a solutions up to $n=10^7$ and with better and better margins. Since everything looked good
at the boundary $b=2$, why look further? But in the row just adjacent to the boundary, there appeared
two exceptional cases, $3+109w$ and $3+121w$ which we call the {\bf Eisenstein ghost twins} as
they appeared out of nothing as the bad guys in the ``matrix" of computations. 
These Eisenstein integers can not be written as a sum of two primes in $Q$. They appear 
topologically isolated, surrounded by "good rows" of the Gaussian integer matrix. So, they appear 
completely unique of this kind: 

\conjecture{
Except for the two Eisenstein ghosts $109w+3,121+3w$, every Eisenstein integer $a+bw$ with $a>1,b>1$ is the 
sum of two Eisenstein primes in $Q$. }

Lets look at the boundary case $z=2+n w$, where we look for two primes $p=1+aw$.
Since $x^2+x+1$ is a cyclotomic polynomial, the Bunyakowski conjecture asking
that there are infinitely many primes of the form $x^2+x+1$ kicks in.
Lets just formulate the boundary case too. It is a special case of enhanced Eisenstein Goldbach 
statement and it isolates the Eisenstein ghost twins. 

\conjecture{
Every integer $n>1$ can be written as a sum $n=a+b$ where $a^2+a+1,b^2+b+1$ are both prime. 
}

\begin{figure}[!htpb]
\scalebox{0.90}{\includegraphics{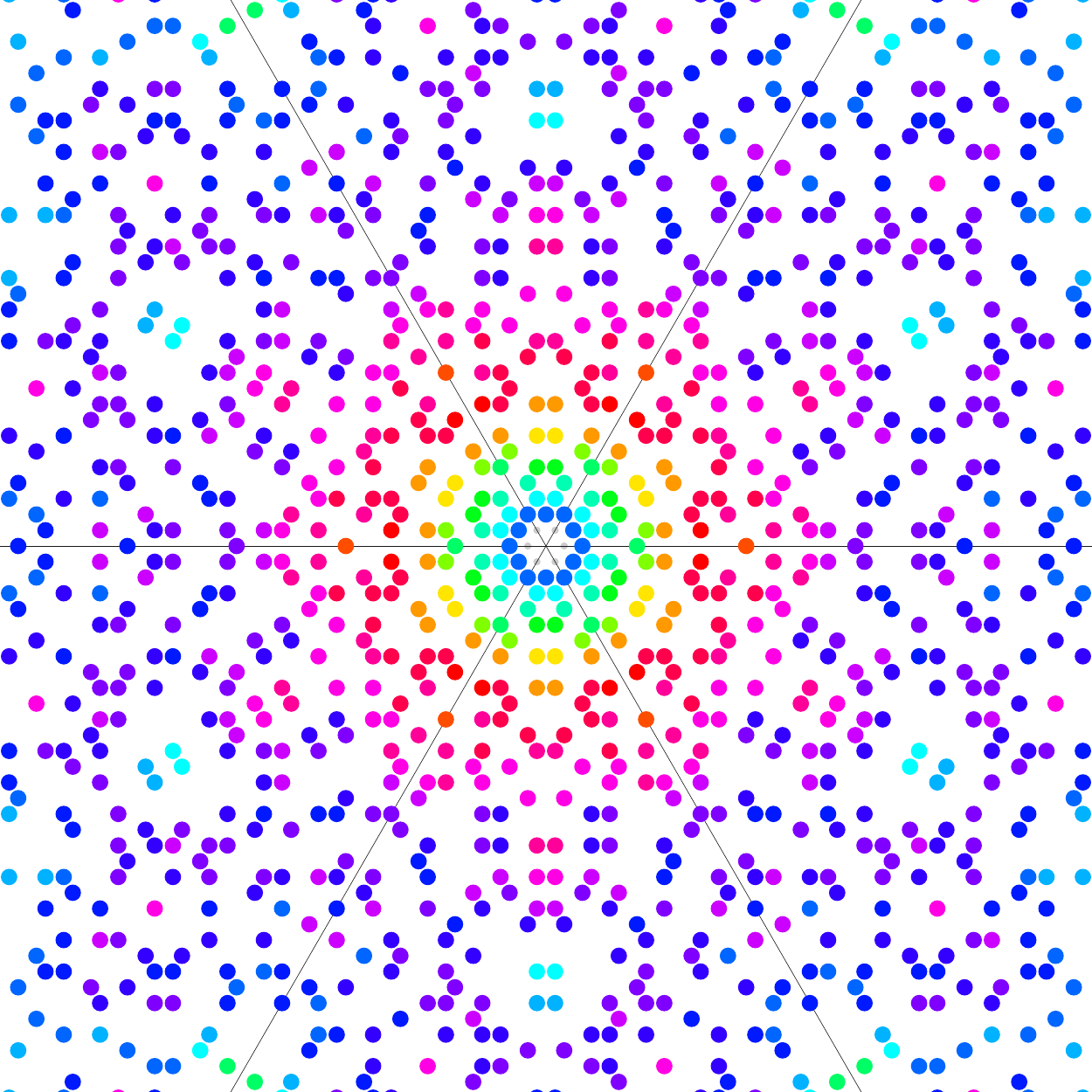}}
\caption{
The Eisenstein primes. 
}
\label{circles}
\end{figure}

Lets state the Mitsui type statement also, as it is so elegant. 

\conjecture{
Every Eisenstein integer can be written as a sum of two Eisenstein primes. 
}

A possible calculus formulation is that the smooth two-periodic function
$f(s,t) = \sum_{p=a+w b} e^{i a t}/a! e^{i b s}/b!$ has the property
that all Fourier coefficients of $g(x,y)=f(x,y)^2$ are positive, if the
sum is over all Eisenstein primes.

\section{Quaternion primes}

Discovered by William Hamilton (1805-1865), the quaternions $z=a+ib+jc+kd=(a,b,c,d)$ have the 
norm $N(z)=a^2+b^2+c^2+d^2$. Their arithmetic is determined by $i^2=j^2=k^2=ijk=-1$,
identities infamously found on October 16, 1843. The set of lattice points with integer values $a,b,c,d$ 
do not form a maximal order. As Adolf Hurwitz (1859-1919) noticed, one needs to include ``half integers". 
Today, integers of the form $(a,b,c,d)$ with integer $a,b,c,d$ are called 
{\bf Lipschitz integers}, and the additional integers of the form 
$(a+1/2,b+1/2,c+1/2,d+1/2)$ with integer $a,b,c,d$ are called the {\bf Hurwitz integers}. 
The union of Lipschitz and Hurwitz integers are then the {\bf Quaternion integers} \cite{Hurwitz1919}. \\

The quaternion integers form a lattice in $R^4$ which is known to be the densest lattice 
sphere packing in $\mathbb{R}^4$, the $D_4$ lattice \cite{KorkinZolotareff} proven by 
Alex Korkin and Yegor Zolotarev (who is also known for the Zolotarev lemma in quadratic reciprocity).
The group of units is the binary tetrahedron group which forms a regular 4-polytop, 
the {\bf 24 cell}, when drawn on the 3-sphere. This is the unit sphere in the discrete $D_4$ lattice.
Despite the lack of commutativity and prime factorization which 
depends on the order in which the primes are arranged, the norm property of division algebra still allows to 
define {\bf quaternion primes} as quaternion integers $p$ which can not be divided by any other 
quaternion integer $q$ different from a unit and having smaller norm $N(q)<N(p)$. 
As noticed also again by Hurwitz, their structure is easier: 
they are just the integers for which $N(q)$ is a rational prime. As the integers, also
the Quaternion primes come now in two classes, the {\bf Lipschitz primes} 
and the {\bf Hurwitz primes} depending on whether they are Lipschitz integers or Hurwitz integers.  \\

\begin{figure}[!htpb]
\scalebox{0.12}{\includegraphics{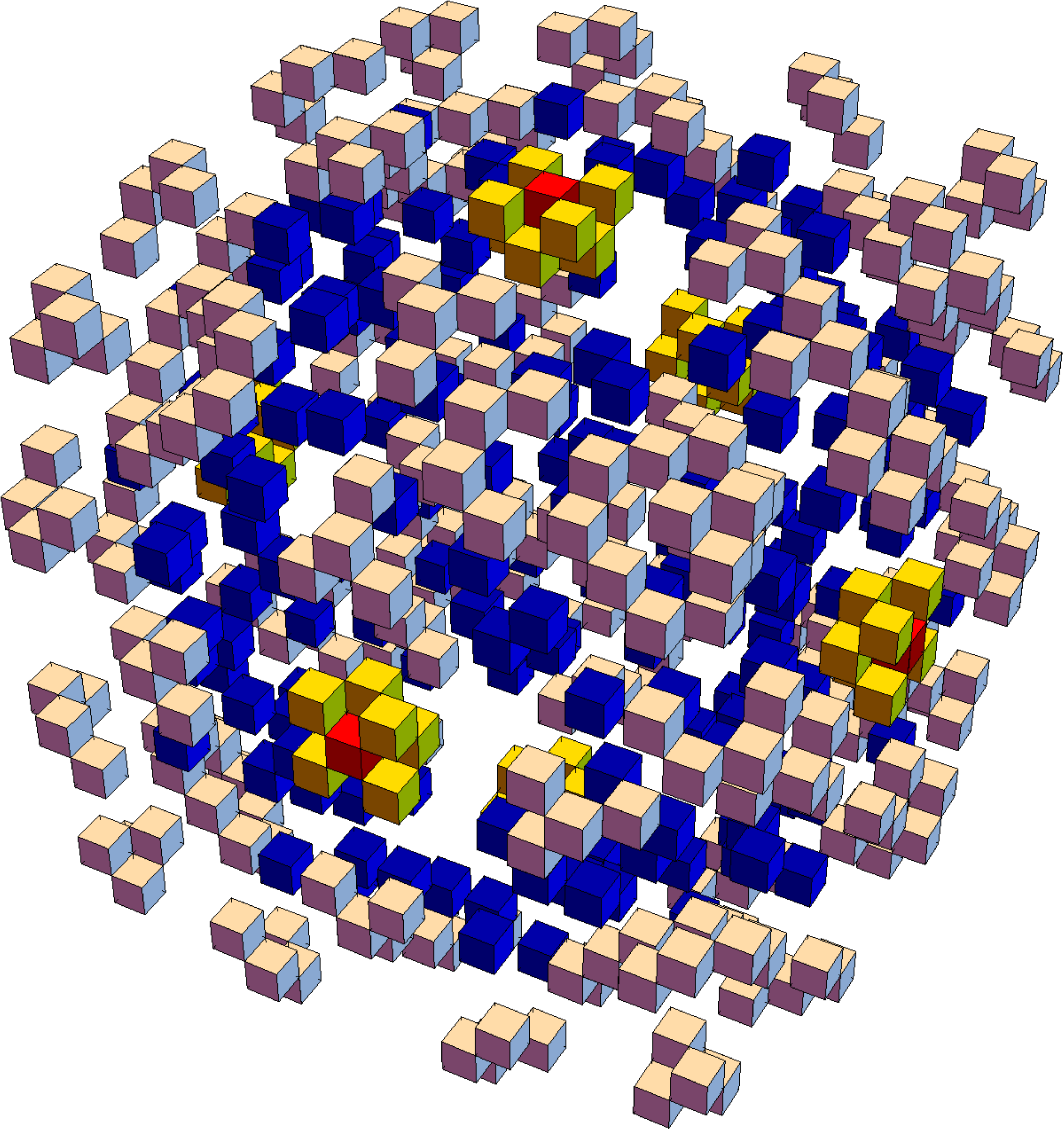}}
\scalebox{0.12}{\includegraphics{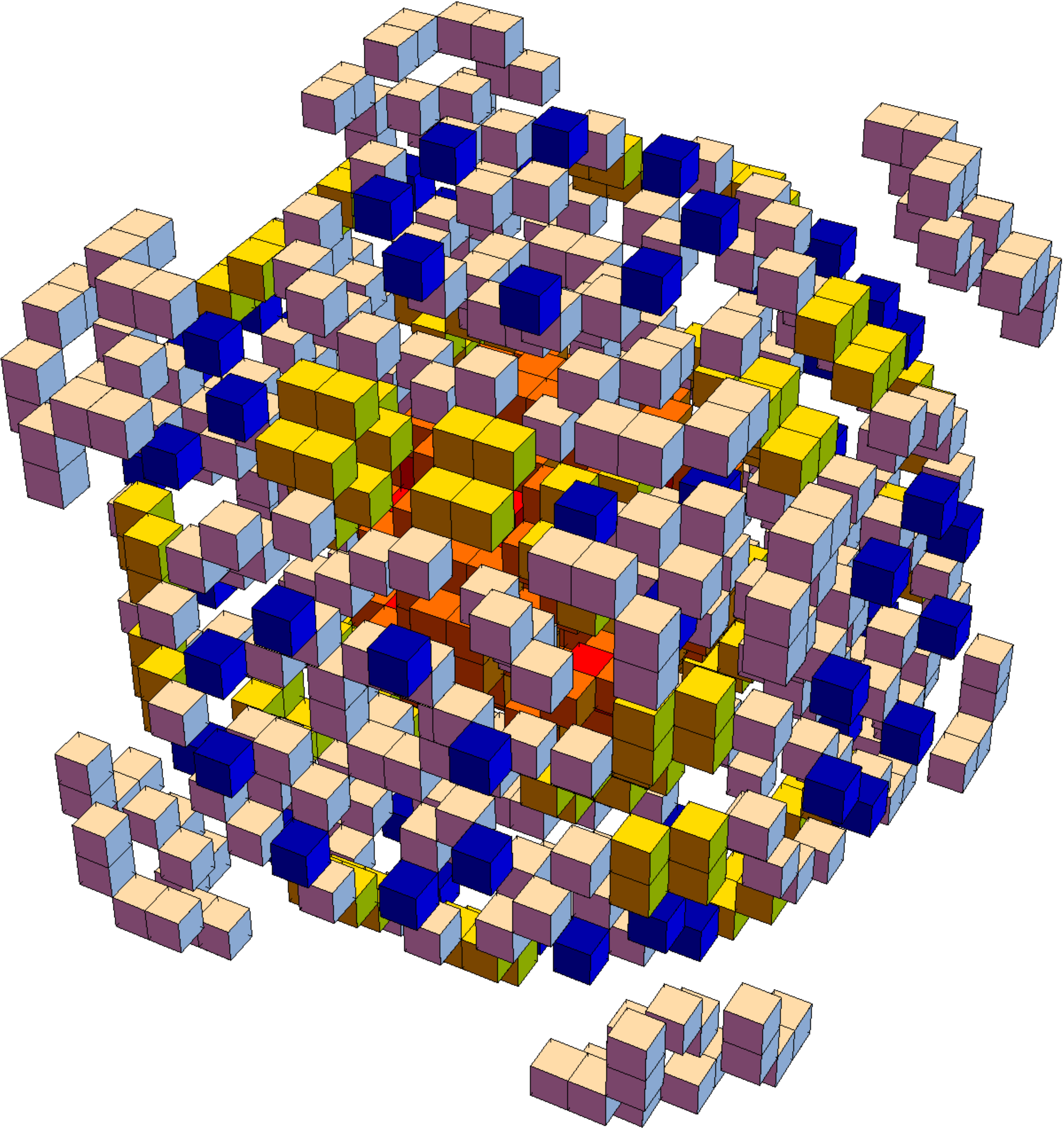}}
\caption{
A slice through all the Lipschitz primes $(a,x,y,z)$ for $a=701$.
They are given by the integer vectors $(a,x,y,z)$ for which $a^2+x^2+y^2+z^2$ is
a rational prime. To the right, we see the Hurwitz primes $(1+2a,1+2x,1+2y,1+2z)/2$ for $a=3001$.
They are given by the integer vectors $(a,x,y,z)+(1,1,1,1)/2$ for which the sum of the
squares is a rational prime.
}
\label{circles}
\end{figure}

\begin{figure}[!htpb]
\scalebox{0.2}{\includegraphics{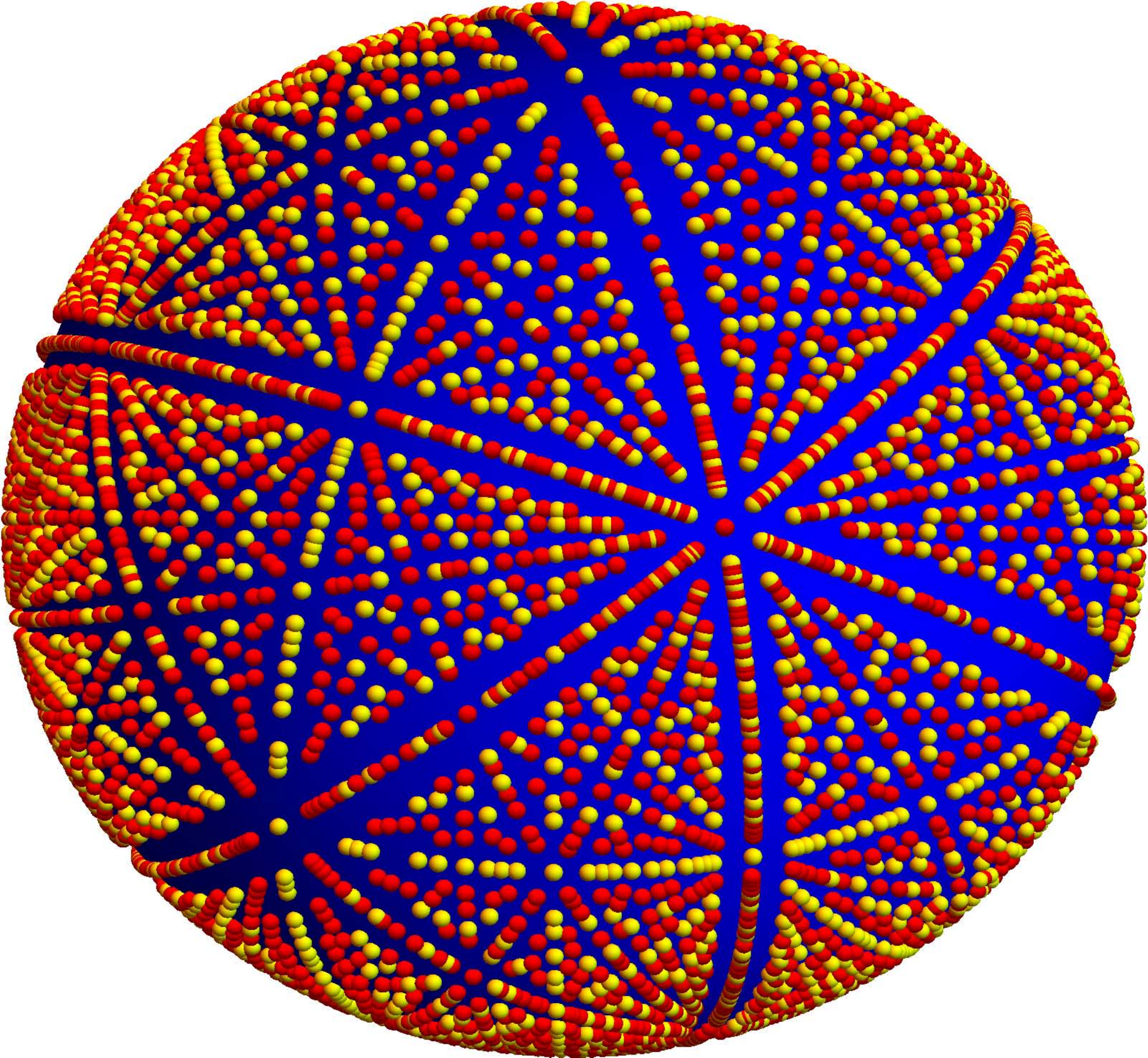}}
\caption{
The {\bf Hurwitz prime sphere} shows the quaternion primes inside the cube $\{ |a|,|b|,|c|,|d| \leq 8\}$ 
projected onto the first three coordinates and then projected onto the unit sphere. 
Since Hurwitz primes can not be on coordinate planes, some grand circles of Lipschitz 
primes appear. Apart from the symmetry imposed by units and conjugation, both Hurwitz 
and Lipschitz prime sets appear to be uniformly distributed on the sphere. }
\label{circles}
\end{figure}

Define the region $Q=\{ (a,b,c,d)  \; | \; a>0,b>0,c>0,d>0 \; \}$. 

\conjecture{Every Lipschitz integer quaternion with entries $>1$ is the sum of two Hurwitz primes in $Q$.}

Again, this could be rewritten analytically as the property that 
$$   f(x,y,z,w) = \sum_{p=a+ib+jc+kd \in Q} x^{2a} y^{2b} z^{2c} w^{2d}  \; , $$ 
summing over all Hurwitz primes in $Q$ has the property that $g=f^2$ has nonzero 
derivatives $g^{(k,l,m,n)}$ for $k,l,m,n>1$. \\

Already for $z=(2,2,2,2)$, there are 14 possibilities:
either $(3,1,1,1)/2 + (1,3,3,3)/2$ (8 cases) or $(1,1,3,3)/2+(3,3,1,1)/2$  (6 cases). 
In the special case, when the Lipschitz integer is $z=(2,2,2,n)$,
the two primes must have up to permutation the form $p = (1,1,1,x)/2,  q=(3,3,3,n-x)/2$,  or then
$p = (1,1,3,x)/2,  q=(3,3,1,2n-x)/2$ for an unknown odd integer $x$.
In the first case, we need simultaneously to have $(3+x^2)/4$  and $(27+(2n-x)^2)/4$ to be rational primes.
In the second case, we need simultaneously to have $(11+x^2)/4$ and $(19+(2n-x)^2)/4$ to be prime. Since $x$
needs to be odd for $p$ to be a Hurwitz prime, we can
write $x=2k+1$. Now $(3+x^2)/4=1+k+k^2$ and $(27+(2n-x)^2)/4 = 7-(n-k)+(n-k)^2$.
In the second case $3+k+k^2$ and $5-(n-k)+(2n-k)^2$. We see:
If the Hurwitz Goldbach conjecture holds, then for any $n$, there exists $k<n$ for which
both $3+k+k^2$ and $7+k+k^2-n-2kn+n^2$ are prime or for which both $1+k+k^2$ and
$5+k+k^2-n-2kn+n^2$ are prime. Now, if Goldbach is true and if there existed only finitely
many primes of the form $3+k+k^2$ and $1+k+k^2$, then for all $m=n-k$ large enough, $7-m+m^2$ and $5-m+m^2$
would always have to be prime. This is obviously not true if $m$ is a multiple of $7$ or $5$. We see that
quaternion Goldbach again implies a special case of the Bunyakovsky conjecture, which like Landau's problem is
likely not so easy to prove:

\remark{
If the quaternion Goldbach conjecture is true, then one of the sequences $1+k+k^2$ or $3+k+k^2$ 
contains infinitely many primes.}

Let us quickly verify the {\bf Bunyakovsky conditions} which enter the {\bf Bunyakovsky conjecture}.
First, $\phi_3(k) = 1+k+k^2$ is already the cyclotomic polynomial and also $k^2+k+3$ satisfies the conditions
of the Bunyakovsky conjecture: the maximal coefficient of the polynomial $f$ is positive, the coefficients have
no common divisor and there is a pair of integers $n,m$ such that $f(n),f(m)$ have no common divisor.
The set of $k$ for which $k^2+k+1$ is prime is the sequence $A002384$ in \cite{A002384}.

\begin{figure}[!htpb]
\scalebox{0.9}{\includegraphics{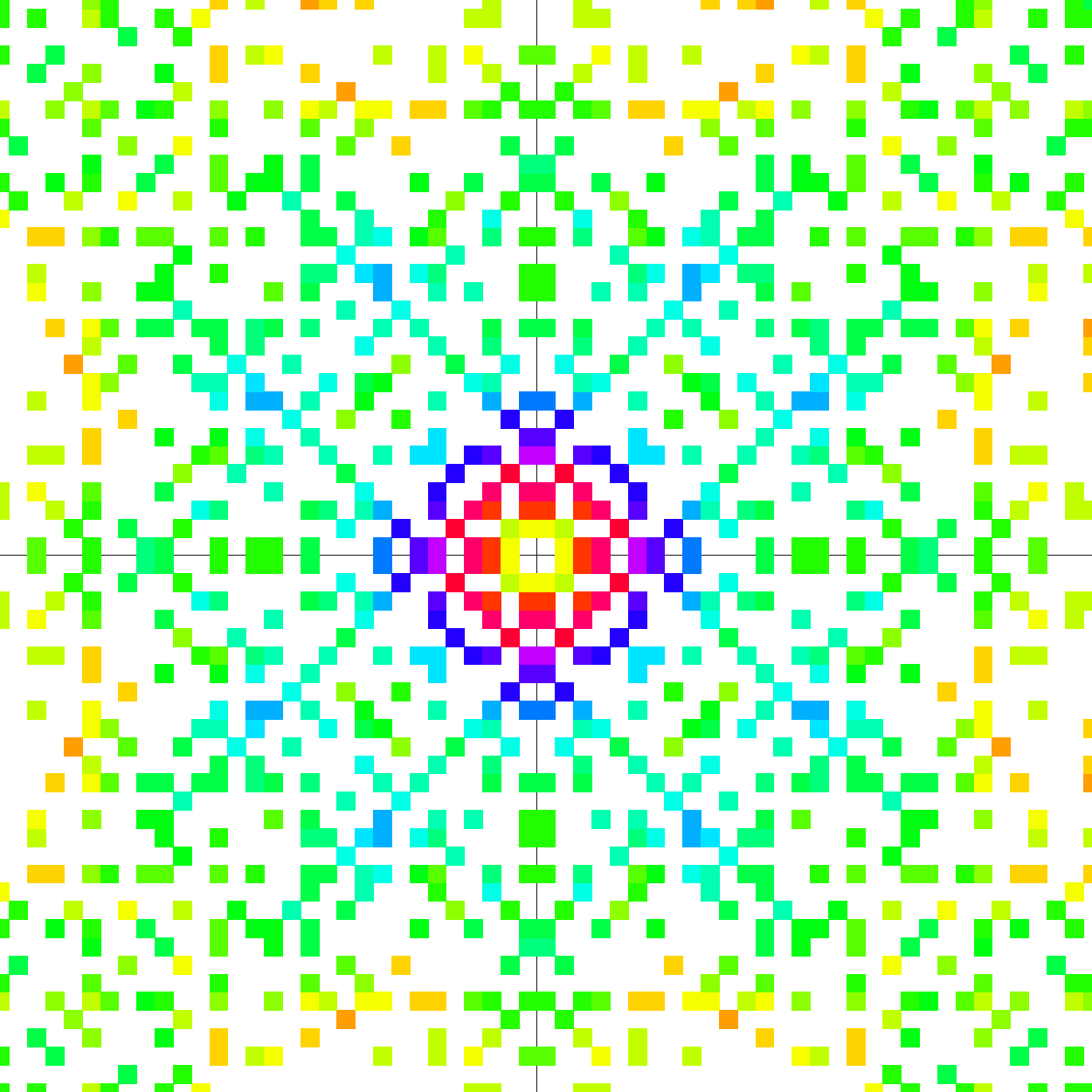}}
\caption{
Hurwitz primes of the form $(1/2,1/2,a+1/2,b+1/2)$. 
}
\end{figure}

We also see experimentally

\conjecture{Every Hurwitz integer quaternion with entries $>2$ is the sum of a Hurwitz and Lipschitz primes in $Q$.}

The Hurwitz integer $(3/2,3/2,3/2,3/2)$ is not the sum of a Hurwitz and Lipschitz prime because the only 
decomposition would be $(1,1,1,1) + (1,1,1,1)/2$ but both are not prime. Together:

\conjecture{Every integer quaternion with entries $>2$ is the sum of two quaternion primes in $Q$.}

\section{Octavian primes}

Besides $\mathbb{R},\mathbb{C}$ and $\mathbb{H}$, there is a fourth normed division algebra
$\mathbb{O}$, the space of {\bf Octonions}. Also called the space of {\bf Cayley numbers} 
or {\bf Hypercomplex numbers}, they were discovered by John Graves and 
not only lack commutativity but have no associativity. 
The members of $\mathbb{O}$ can either be written as a linear combination of a basis
$1,i,j,k,l,m,n,o$ or then, according to a suggestion of Arthur Cayley and Leonard Dickson, 
as pairs $(z,w)$ of quaternions, defining $(z,w) \cdot (u,v) = (zu-v^*w,vz+wu^*)$. 
The algebra is no more associative. \\

In order to do number theory, one has to specify what the {\bf integers} are in $\mathbb{O}$.
There are now three classes of integers, the {\bf Gravesian integers} $(a,b,c,d,e,f,g)$,
the {\bf Kleinian integers} $(a,b,c,d,e,f,g) + (1,1,1,1,1,1,1,1)/2$ as well as the
{\bf Kirmse integers} which includes elements for which $4$ of the entries are half integers.
There are 7 maximal orders which Johannes Kirmse classified in 1925 \cite{Kirmse25,Coxeter46}.
They are all equivalent and produce
the class of {\bf Octonion integers} or {\bf Cayley integers} or then more catchy, 
the {\bf Octavian integers}. See \cite{ConwaySmith,Baez2002}. 
The condition $N(z w) = N(z) N(w)$ which assures that the algebra a {\bf division algebra}, 
is also called the {\bf Degen eight square identity} named after Carl Ferdinand Degen.
The identity allows to define Octavian primes as the set of Octavian integers for
which the sum of the squares is a rational square.  

\begin{figure}[!htpb]
\scalebox{0.2}{\includegraphics{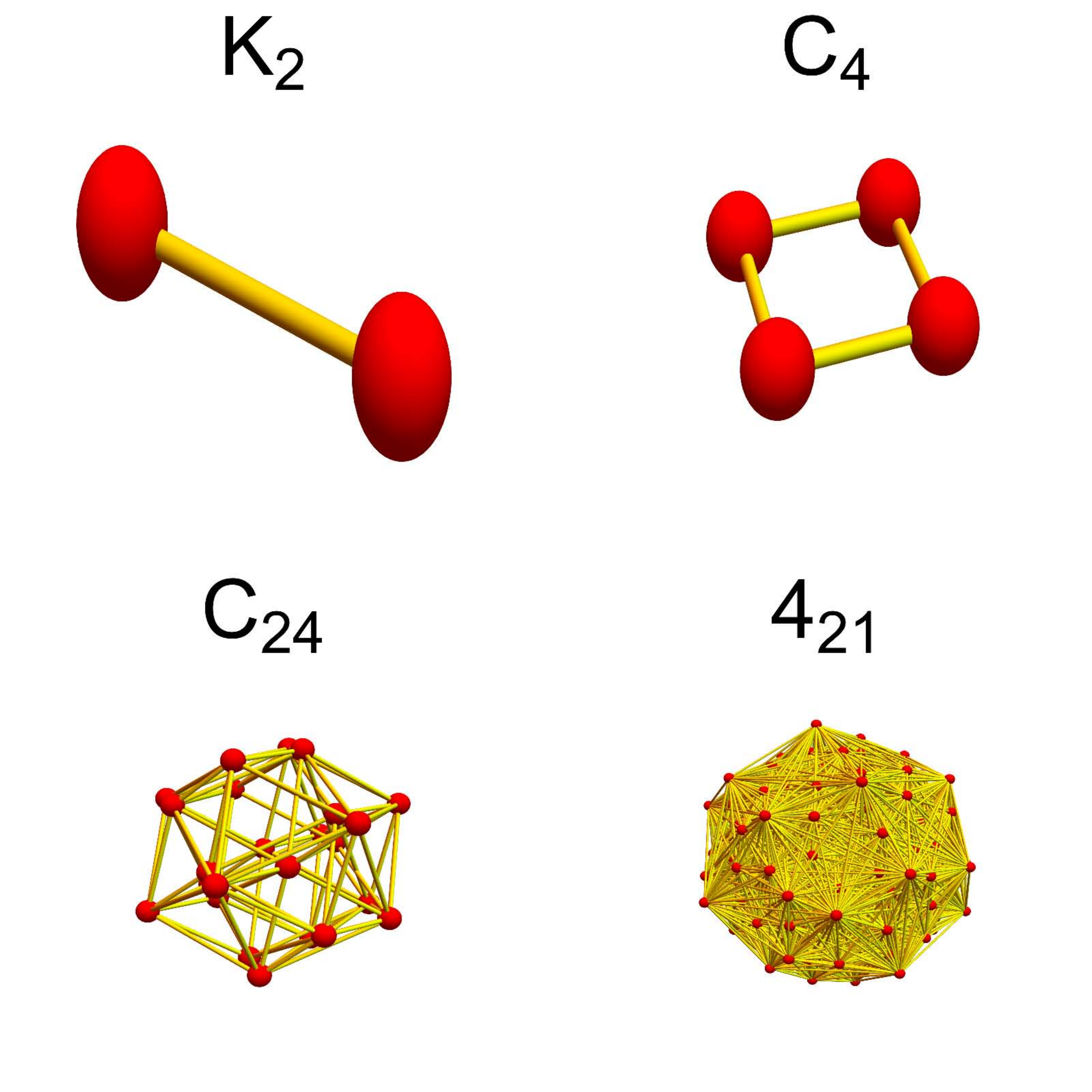}}
\caption{
The units in the four normed division algebras
$\mathbb{R},\mathbb{C},\mathbb{H},\mathbb{O}$: We see the $K_2$ for $\{ \pm 1\} \subset \mathbb{R}$,
$C_4 = \{1,i,-1,-i\} \subset \mathbb{C}$, the 24 cell in $\mathbb{H}$ and the Gosset polytop 
in $\mathbb{Q}$ generating the $E_8$ lattice. The all are known to generate the densest sphere packings,
except for $\mathbb{C}$, where one has to look at the Eisenstein integers instead.
\label{division}
}
\end{figure}

Primes can now be called {\bf Gravesian}, {\bf Kleinian} or ${\bf Kirmse}$
primes. All three distinct classes together form the {\bf Octavian primes}.
Their {\bf units}, the unit norm integral Octonions, are remarkable.
They they do not form a multiplicative {\bf group}, as multiplication is not associative,
they form a finite {\bf loop} of 240 elements. A loop is an algebraic structure 
which is more primitive and so general than a group in which one does not insist on associativity. 
The units form what is one calls a {\bf Moufang loop} \cite{Moufang33} named after German mathematician 
{\bf Ruth Moufang} (1905-1977) who was a student of Max Dehn (1878-1952), who
is known for the Dehn-Sommerville relations or Dehn twists. There is a smaller loop of 16 unit 
octonions containing the Gravesian integers like $(\pm 1,0,0,0,0,0,0,0,0)$. 
The units placed in the unit sphere of $R^8$ form the {\bf Gosset polytope} $4_{21}$ which 
was discovered by Thorold Gosset (1869-1962), who as a lawyer without much clients
amused himself as an amateur mathematician and also helped proofreading
\cite{coxeter} which has a note on Gosset on page 164. The vertices of $4_{21}$ are the 
roots of the exceptional Lie algebra $E_8$ belonging to the 248 dimensional Lie group 
$E_8$. As the dimension of the maximal torus is $8$, the root system lives in $R^8$.
One can write points on the sphere of radius $2$ taking vertices $(a,b,0,0,0,0,0,0)$
with $a,b \in \{-1,1\}$ or $(a,b,c,d,e,f,g,h)$ with with entries in $\{-1/2,1/2\}$
summing up to an even number. This lattice $E_8$ has just recently been verified by 
Maryna Viazovska to be the densest sphere packing in $R^8$ \cite{Viazovska}. \\

In the Octavian case, it appears that it is not obvious how to come up with a conjecture which 
both convinces and is justifiably difficult. Brute force searches are difficult as the volume 
of a box of size $r$ grows like $r^8$. Here is a first attempt of get to a conjecture. 
We formulate as a question since our experiments did not get far yet, nor do we have an 
idea how difficult the statement could be. Anyway, lets denote by $Q$ again the set of all 
integer octonions, for which all coordinates are positive. 

\question{
Is there a constant $K$ such that every Octavian integer $Z$ with
coordinates $\geq K$ is a sum of two Octavian primes $P,Q$ in $Q$ ? }

While $(1,1,1,1,1,1,1,1)$ as a sum of $p+p$ with Kleinian prime $p=(1,1,1,1,1,1,1,1)/2$ of norm $2$,
already the Kleinian $(1,1,1,1,1,1,1,2)$ can not be written as a sum of two primes
and the Kirmse integer $(2,2,3,3,3/2,3/2,3/2,3/2)$ can not be written as as sum of two primes
the constant has to be at least $2$. Lets look at some small examples using the notation $\overline{a}$
for a block of $4$ numbers $a$. The Kleinian $(\overline{3},\overline{3})/2$ is the 
sum of the two Kirmse primes $(\overline{1/2},\overline{1})$ + $(\overline{1},\overline{1/2})$ 
of norm $5$. The Kirmse integer $(\overline{1},\overline{3/2})$ is the sum 
$(\overline{1/2},\overline{1/2}) + (\overline{1/2},\overline{1})$ of a Kleinian prime of norm $2$ and
a Krimse prime of norm $5$. The integer $(\overline{2},\overline{2})$ can be written as a sum
$(\overline{1/2},\overline{1})$+$(\overline{3/2},\overline{1})$ of two Kirmse primes with 
norm $5$ and $13$. And $(\overline{2},\overline{3/2})$ can be written as a sum
$(\overline{1/2},\overline{1/2})$+ $(\overline{3/2},\overline{1})$ of a Kleinian and Krimse
prime of norm $2$ and $13$. The Kleinian integer $n=(\overline{3},\overline{3})$
already can be written in 64 ways as a sum of two positive Gravesian primes, like $n=p+q$ with 
$p=(\overline{2},2,2,2,1)$ and $q=(\overline{1},2,2,2,1)$ and then for example 6 more ways as a sum of 
Kirmse primes like $p=(1,1,2,2,\overline{3/2})$ and $q=(2,2,1,1,\overline{3/2})$. 
For $(\overline{3},3,3,3,4)$ we have a decomposition like 
$(1,2,2,1,3/2,3/2,3/2,5/2)$,$(2,1,1,2,\overline{3/2})$. \\

\begin{figure}[!htpb]
\scalebox{0.9}{\includegraphics{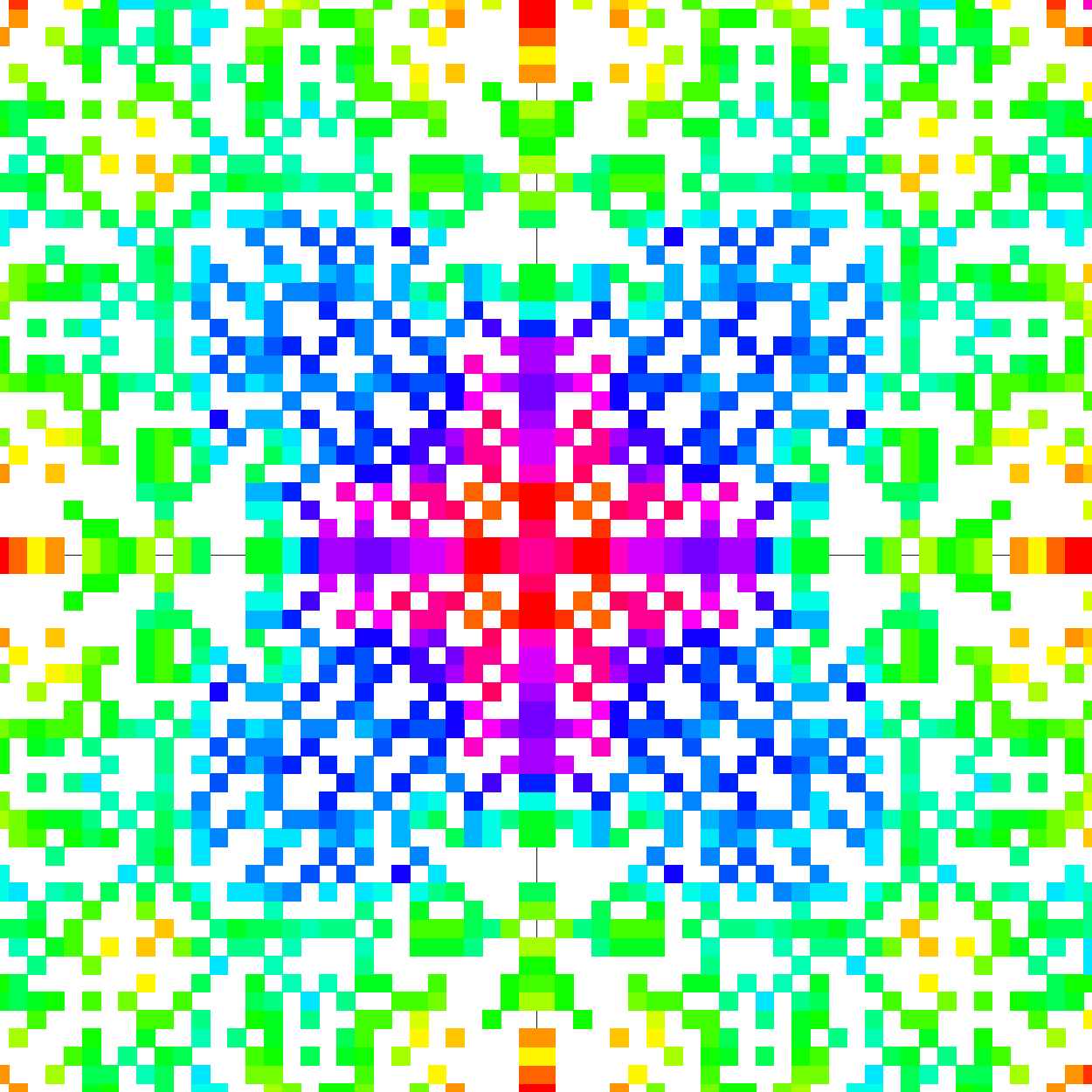}}
\caption{
Kirmse primes $(1,1,1,1,3/2,3/2,a+1/2,b+1/2)$.
They are given by all pairs $(a,b)$ for which $17+a^2+b^2+a+b$ is 
a  rational prime.
\label{kirmse}
}
\end{figure}

To simplify, one can take the Cayley-Dickson point of view and see an Octonion
as a pair $(z,w)$ of quaternions and $(z,w) \cdot (u,v) = (zu-v^*w,vz+wu^*)$, where $z^*=(a,-b,-c,-d)$ 
is the conjugate quaternion of $z=(a,b,c,d)$. A subclass of Octonion prime are now pairs $(z,w)$ of
quaternion integers for which $N(z)+N(w)$ is prime. These could already be enough to work with and get the sum. 
We then don't even have to touch the multiplication of Octonions as the prime property is visible from the 
norm. One certainly has to distinguish cases. In the simplest case, within Gravesian integers
if $(z,w)$, we want to find a Gravesian prime $(x,y)= (x,y)$
such that both $p=\sum_i x_i^2 +y_i^2$ and $q=\sum_i (z_i-x_i)^2 + (w_i-y_i)^2$ are prime. 
One can also just look at two dimensional slices in the 8-dimensional space and notice that the
primes are quite ``dense". For Figure~(\ref{kirmse}) for example, we show at all pairs $(a,b)$ for 
which $a^2+b^2+a + b + 5$ is prime. Note that they appear dense but proving even that some exist
in each row is a Landau-Bunyakovsky type problem. \\

Maybe there is a trick to write down the primes $p,q$ directly using the Lagrange 4 square theorem
so that they add up to $n$. Maybe also that in the higher dimensional situation the original 
{\bf Schnirelman density} approach works. This density $\alpha$ is a quantity of a set $A$ such that
for $X=[-M,M]^8$ satisfies $|A \cap X(M)| \geq \alpha |X(M)| = \alpha (2M+1)^8$ for every $M>1$. 
Schnirelman showed in the one dimensional case that the density of the set sum $P+P$ of primes $P$
is positive. In the higher dimensional case, a similar argument could be used. 
In any way, the intuitive argument given by Hardy and Littlewood for primes makes things
plausible: in all cases, the probability to find a prime in the box $Q(n)$ is $1/\log(n)$. 
To hit two primes has probability $1/\log(n)^2$, but we have $|Q(n)|$ chances. 
The chance to miss in higher dimensions shrinks even faster with $n$ in higher dimensions. \\

But we don't even have scratched the surface with experimental investigations in 
the case of Octavian integers. \\

{\bf Acknowledgements}: Some computations were run on the Odyssey
cluster supported by the FAS Division of Science, Research Computing Group
at Harvard University.

\section*{Appendix: Goldbach comets and Ghosts in the matrix }

For Gaussian Goldbach, things look comfortable. It looks even
possible that one can chose for every row a small set of rows.

\begin{figure}[!htpb]
\scalebox{0.4}{\includegraphics{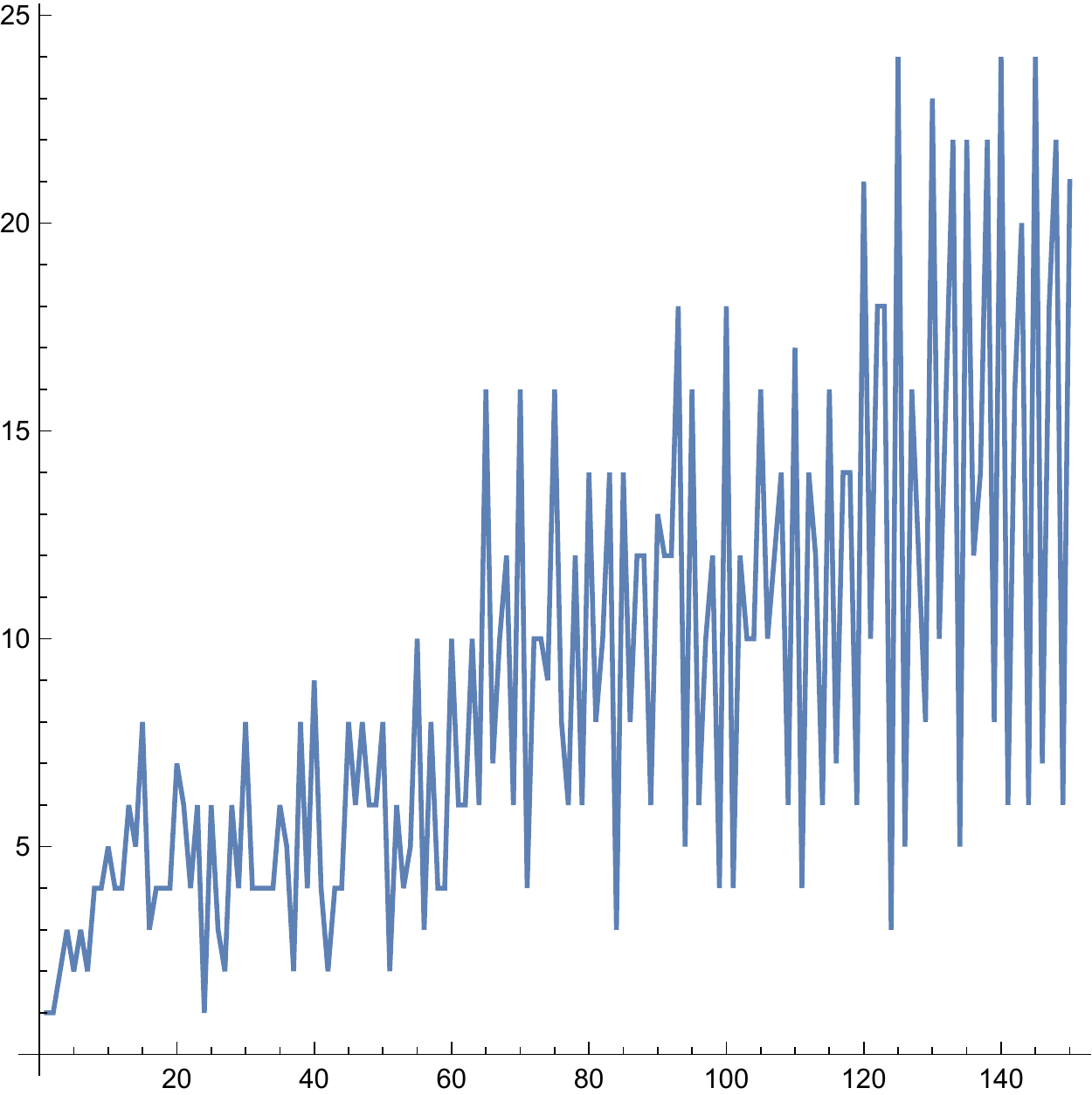}}
\scalebox{0.4}{\includegraphics{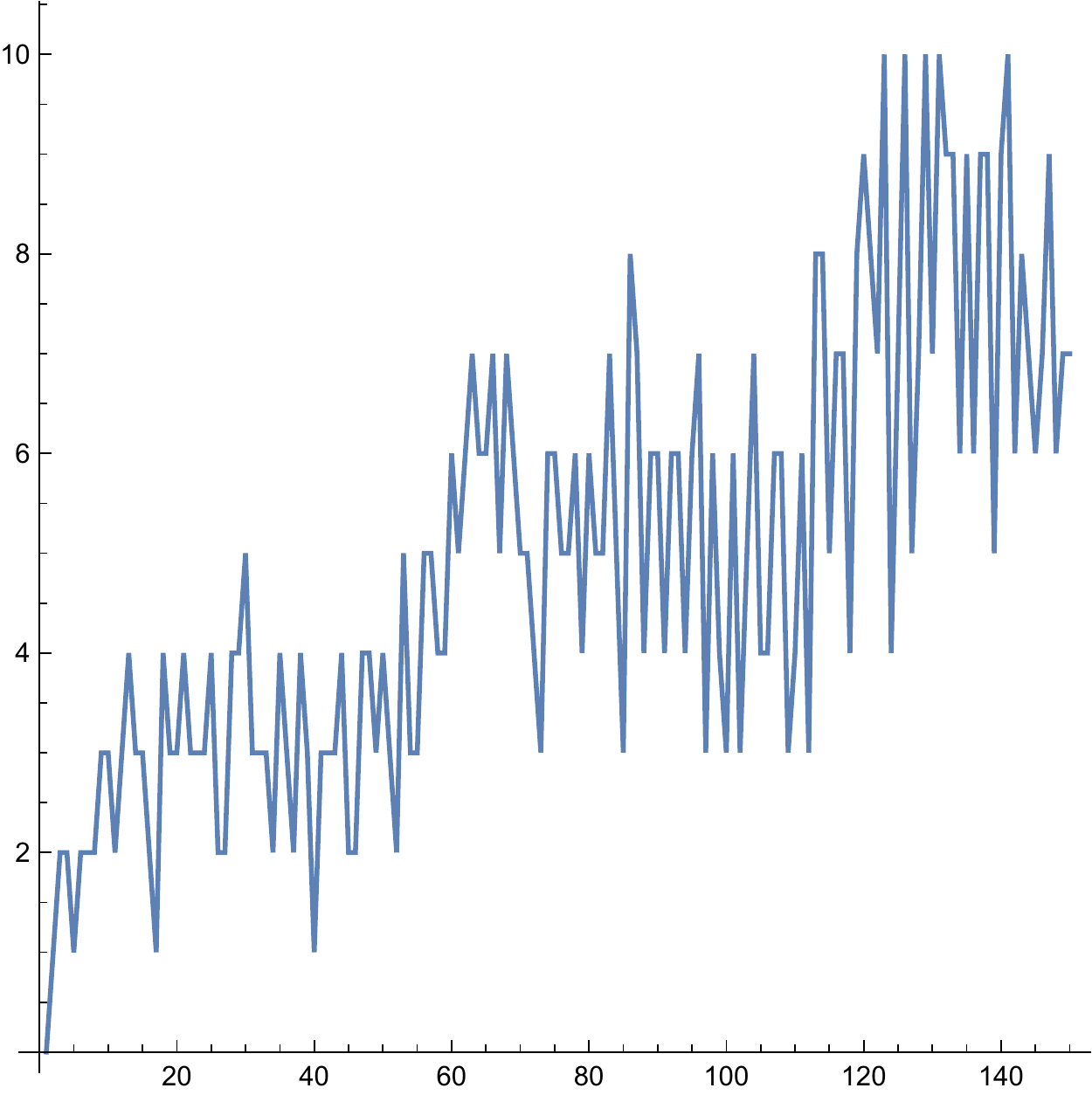}}
\scalebox{0.4}{\includegraphics{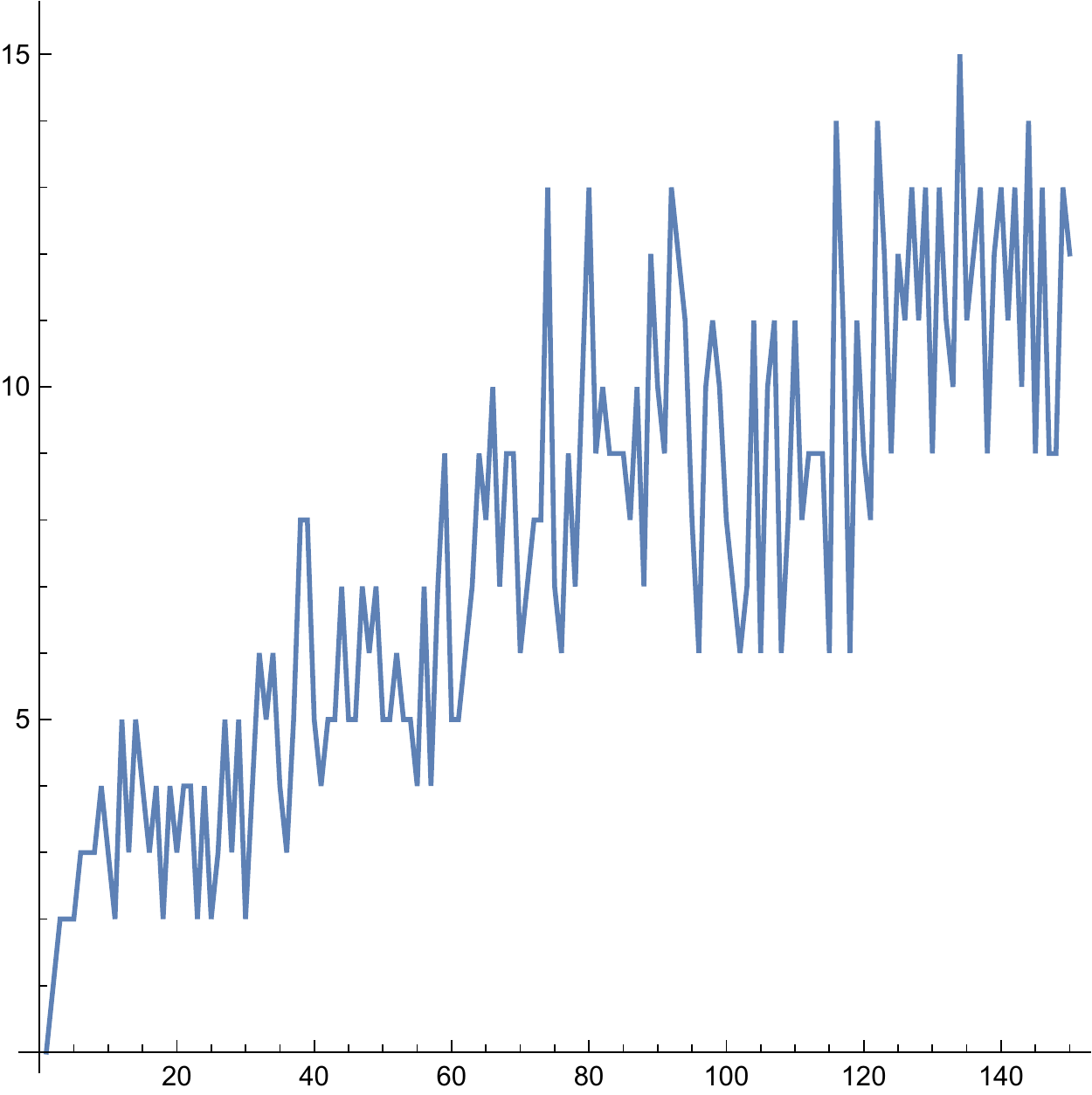}}
\scalebox{0.4}{\includegraphics{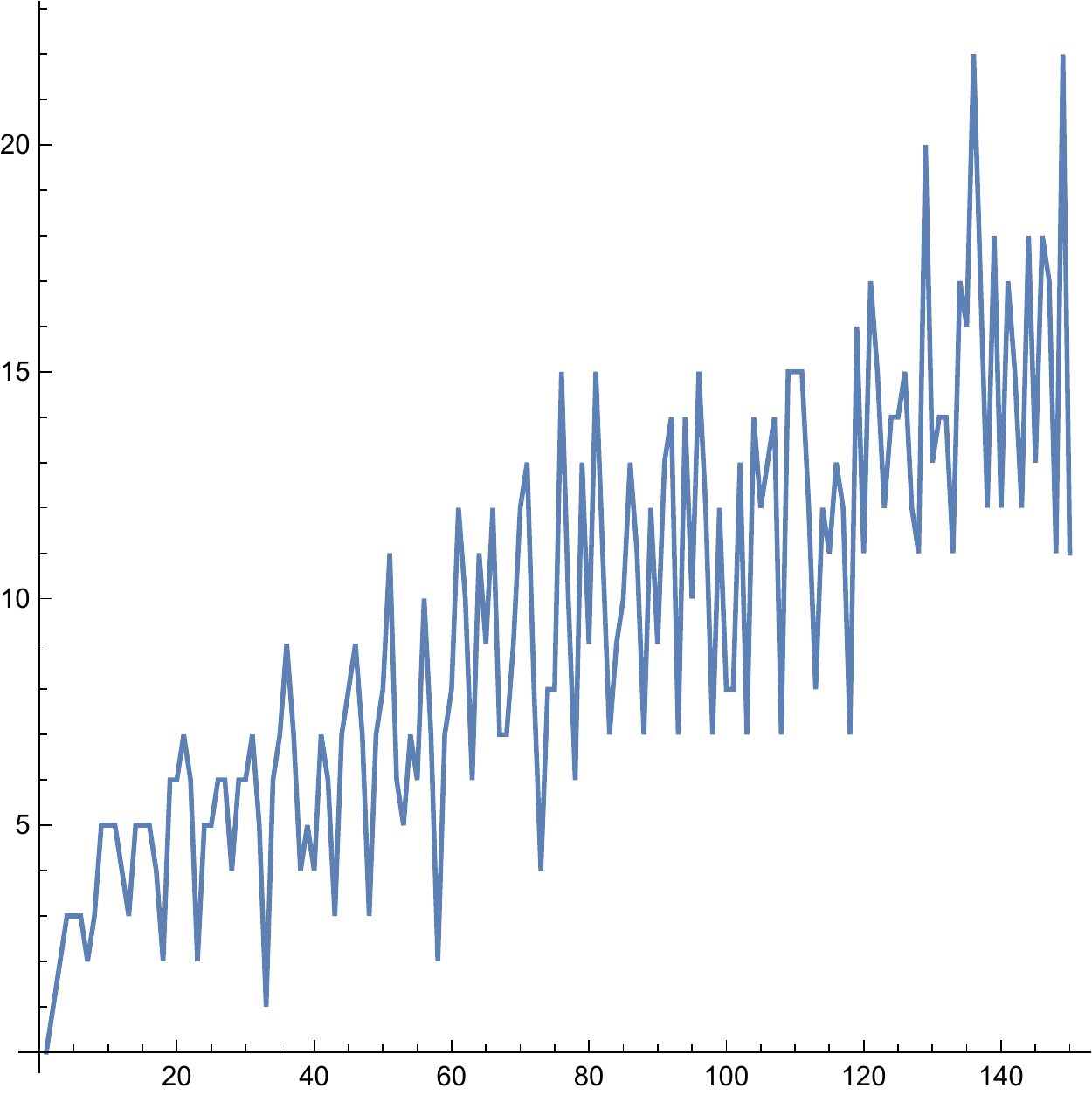}}
\scalebox{0.4}{\includegraphics{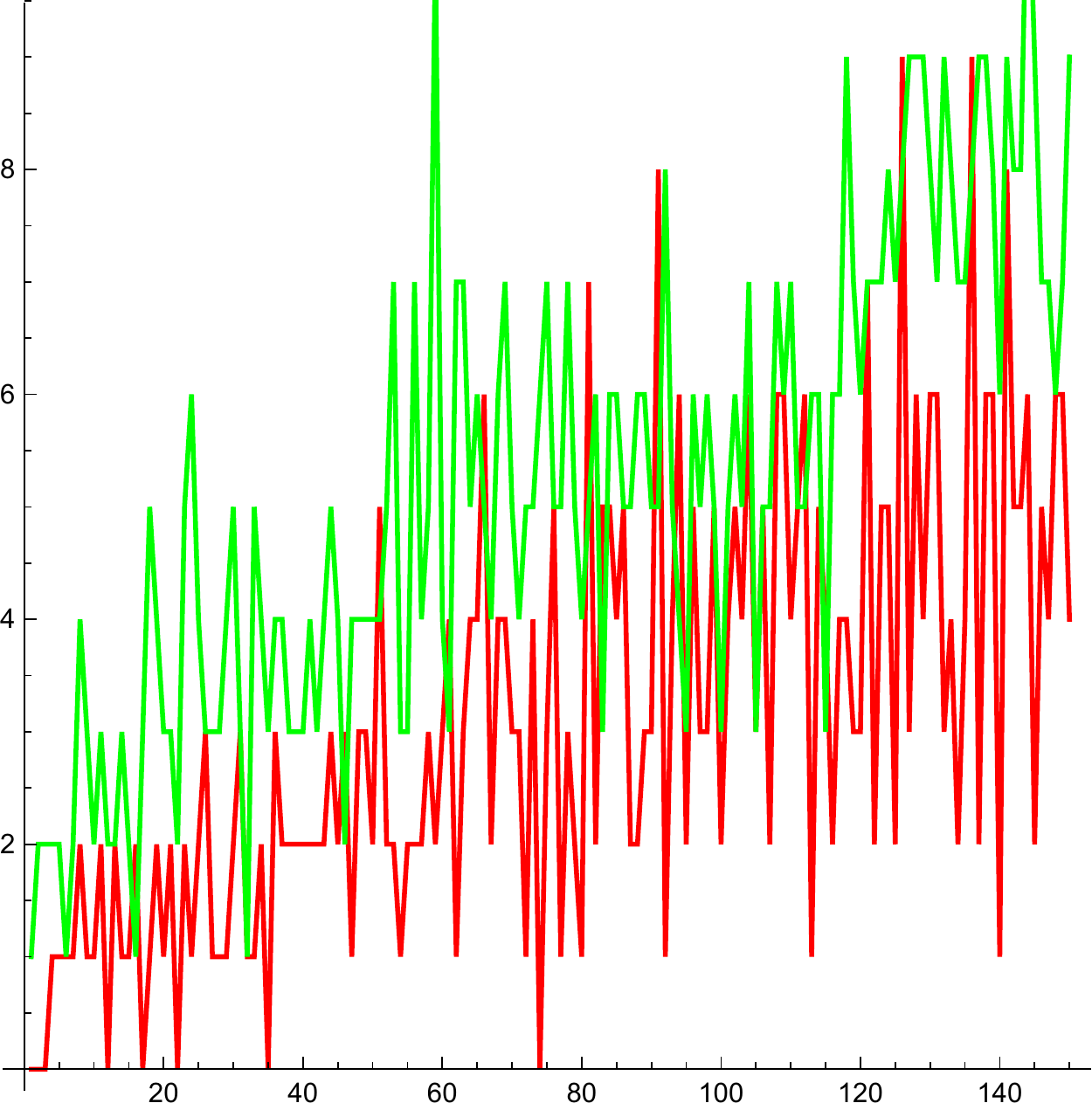}}
\scalebox{0.4}{\includegraphics{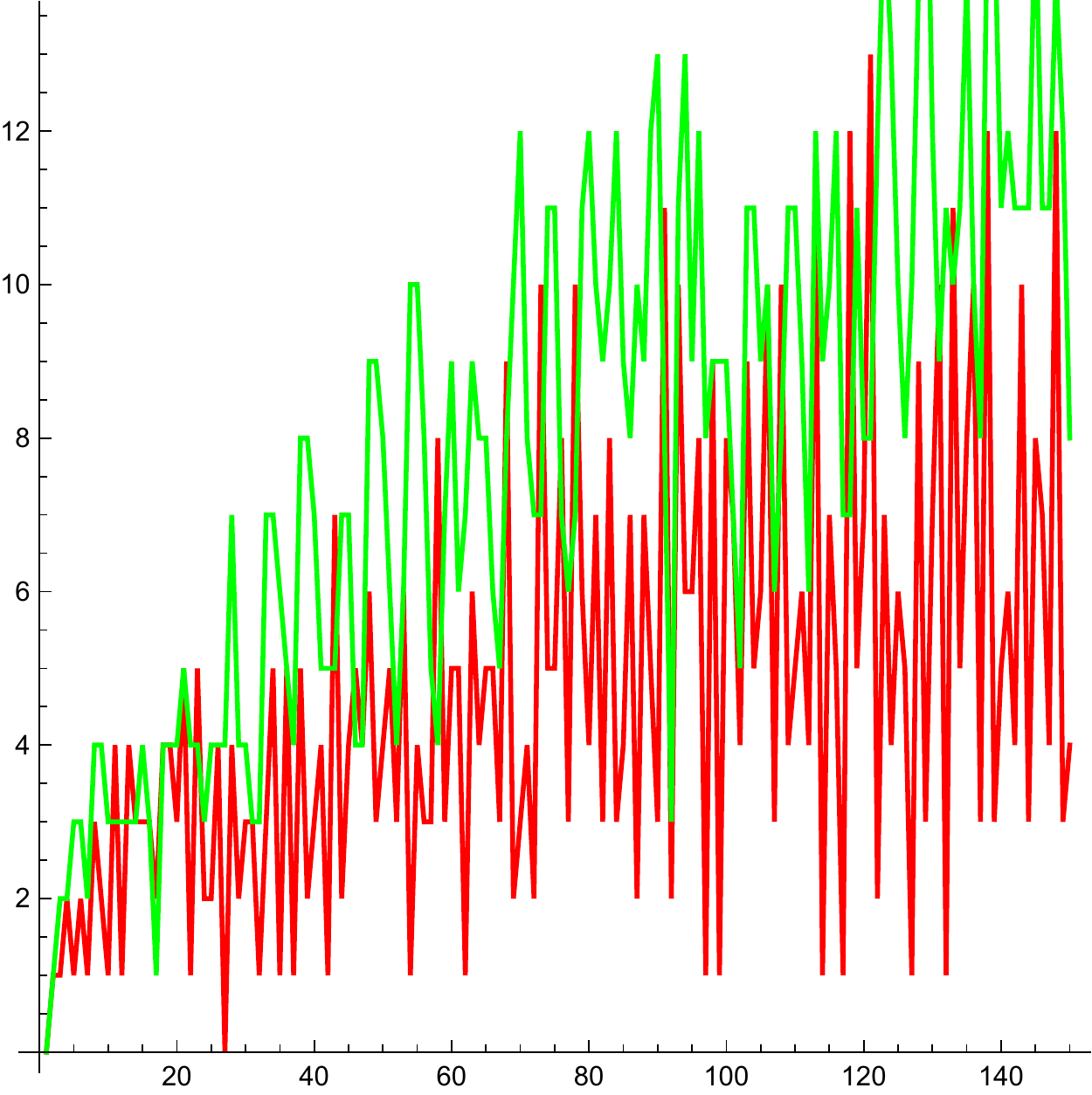}}
\caption{
The number of times, a Gaussian integer $a+ i b$ can
be written as a sum of Gaussian primes. We see cases
$b=2,b=3,b=4,b=5,b=6,b=7$. In the last two cases, one
can not force a prime from the first row.
}
\end{figure}

We next now pictures illustrating
{\bf Eisenstein ghost twins} with "ghost" added, not to
confuse with Eisenstein prime twins, which are neighboring
Eisenstein primes.

\begin{figure}[!htpb]
\scalebox{0.4}{\includegraphics{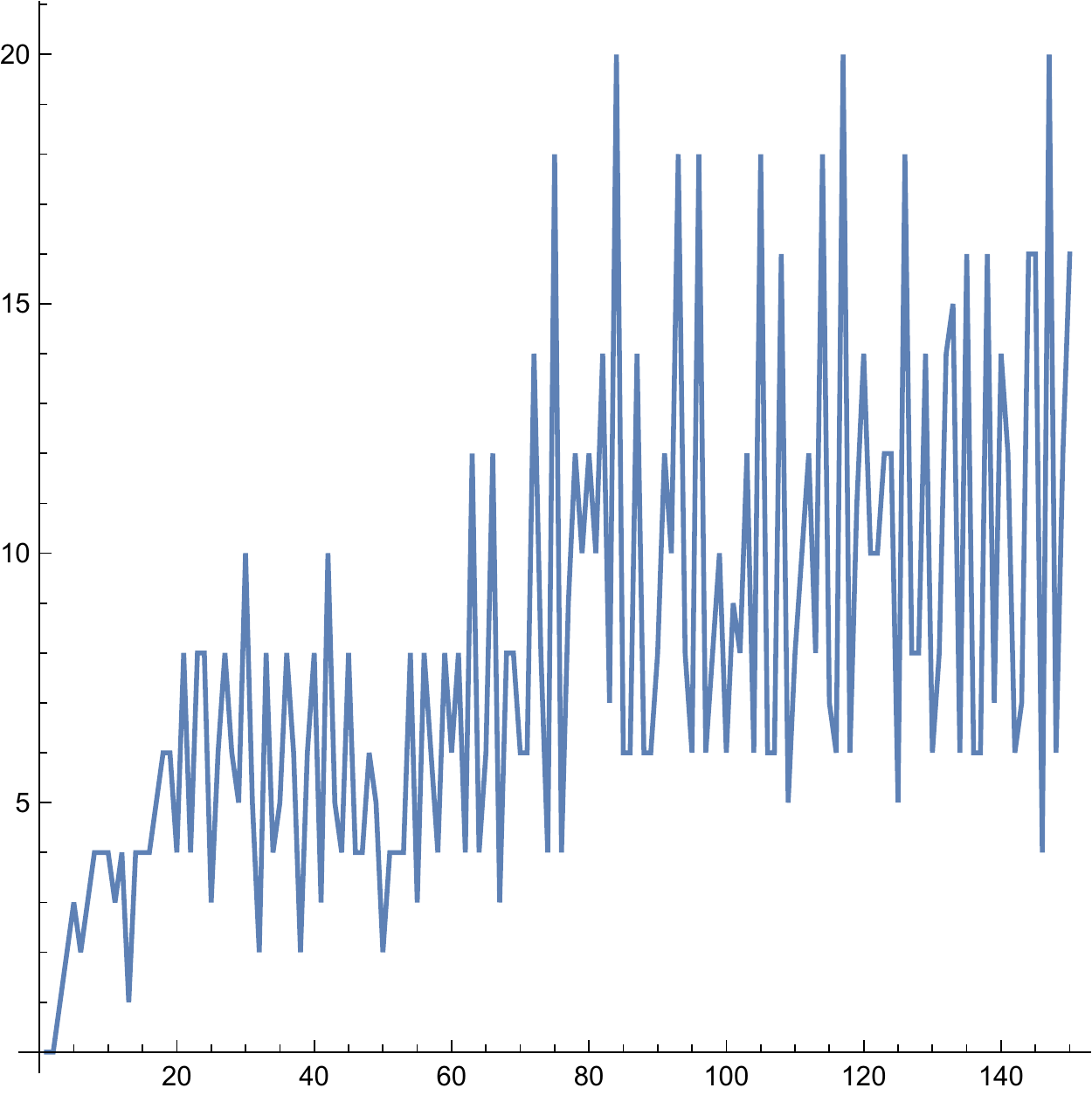}}
\scalebox{0.4}{\includegraphics{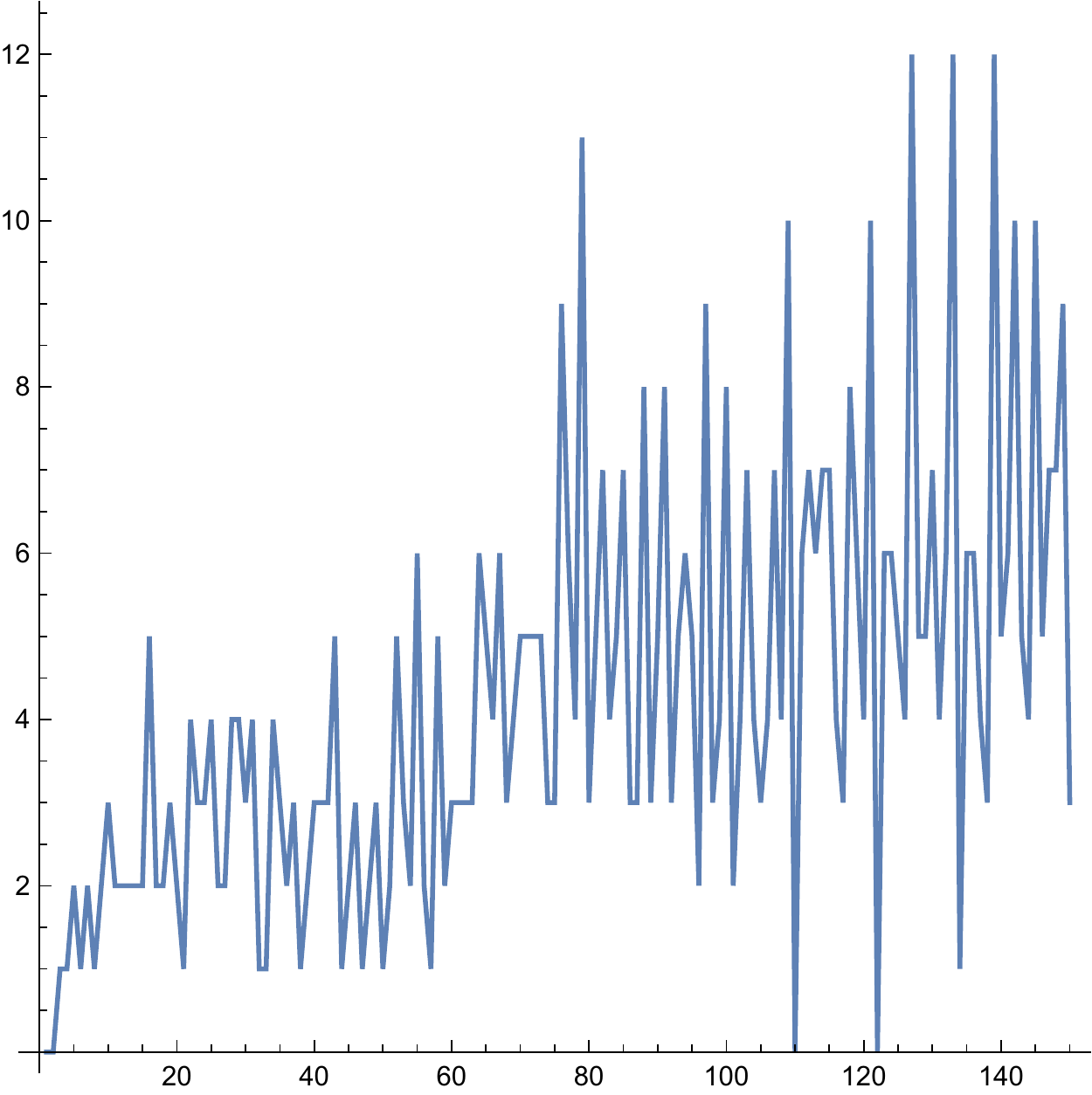}}
\scalebox{0.4}{\includegraphics{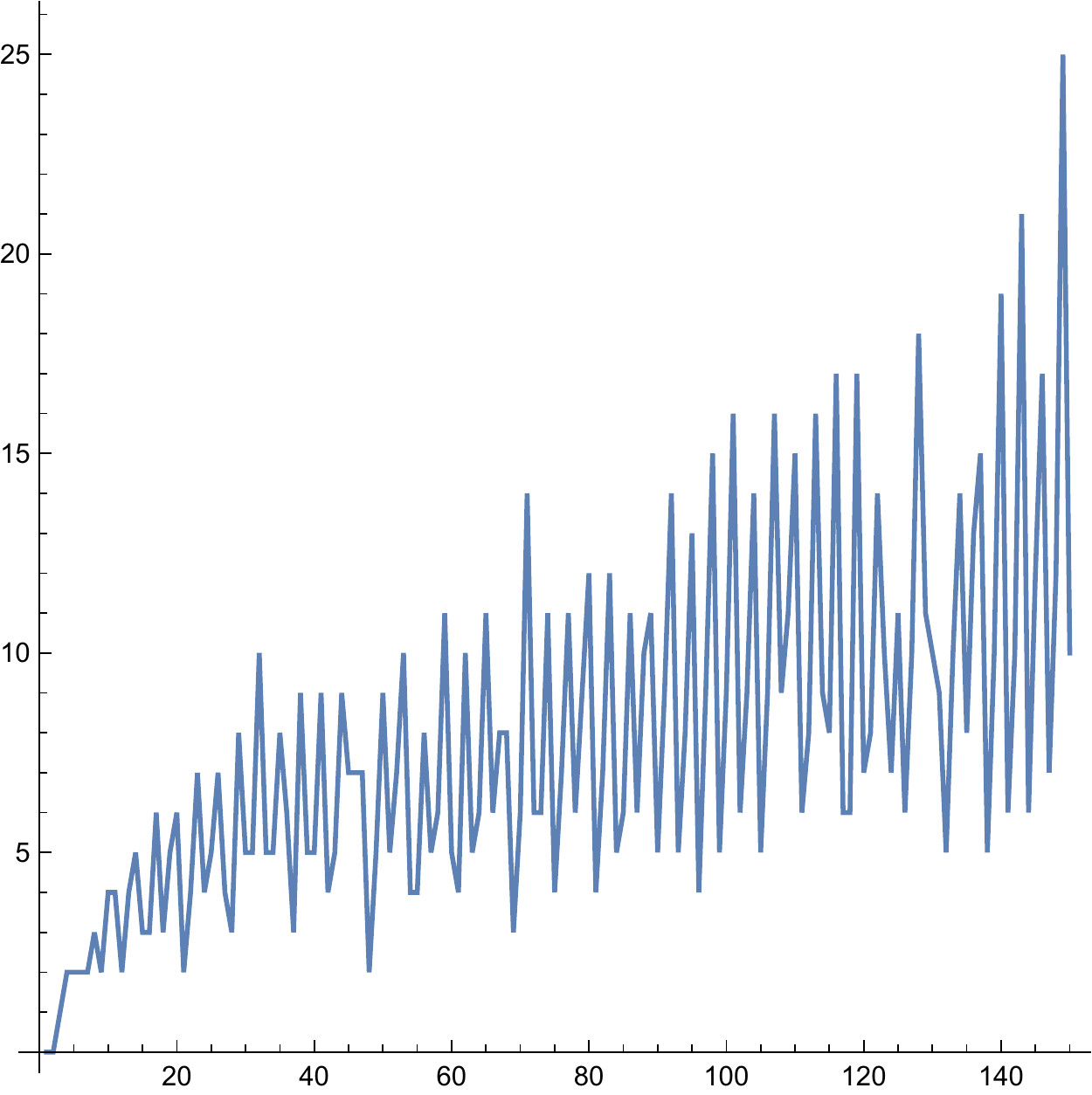}}
\scalebox{0.4}{\includegraphics{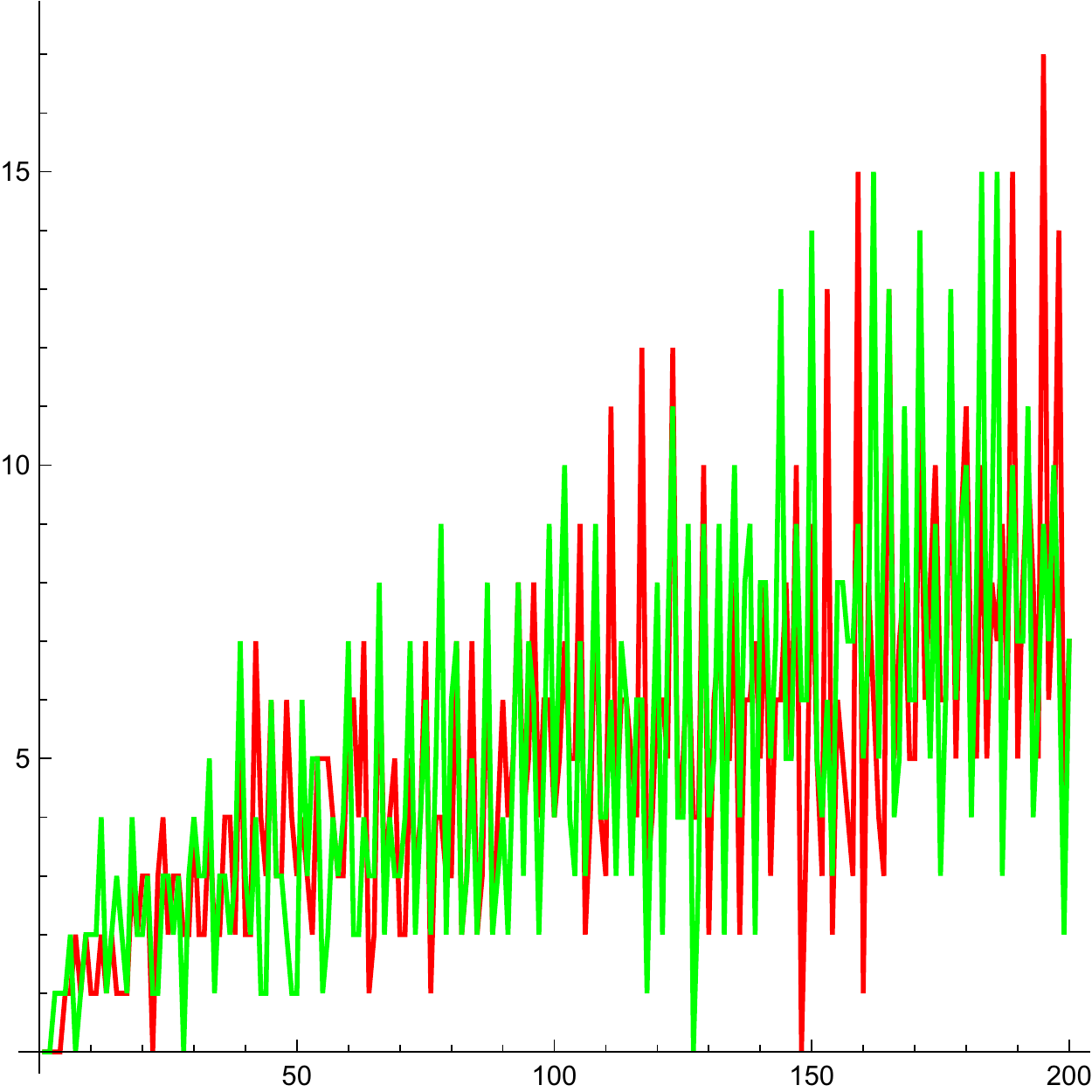}}
\scalebox{0.4}{\includegraphics{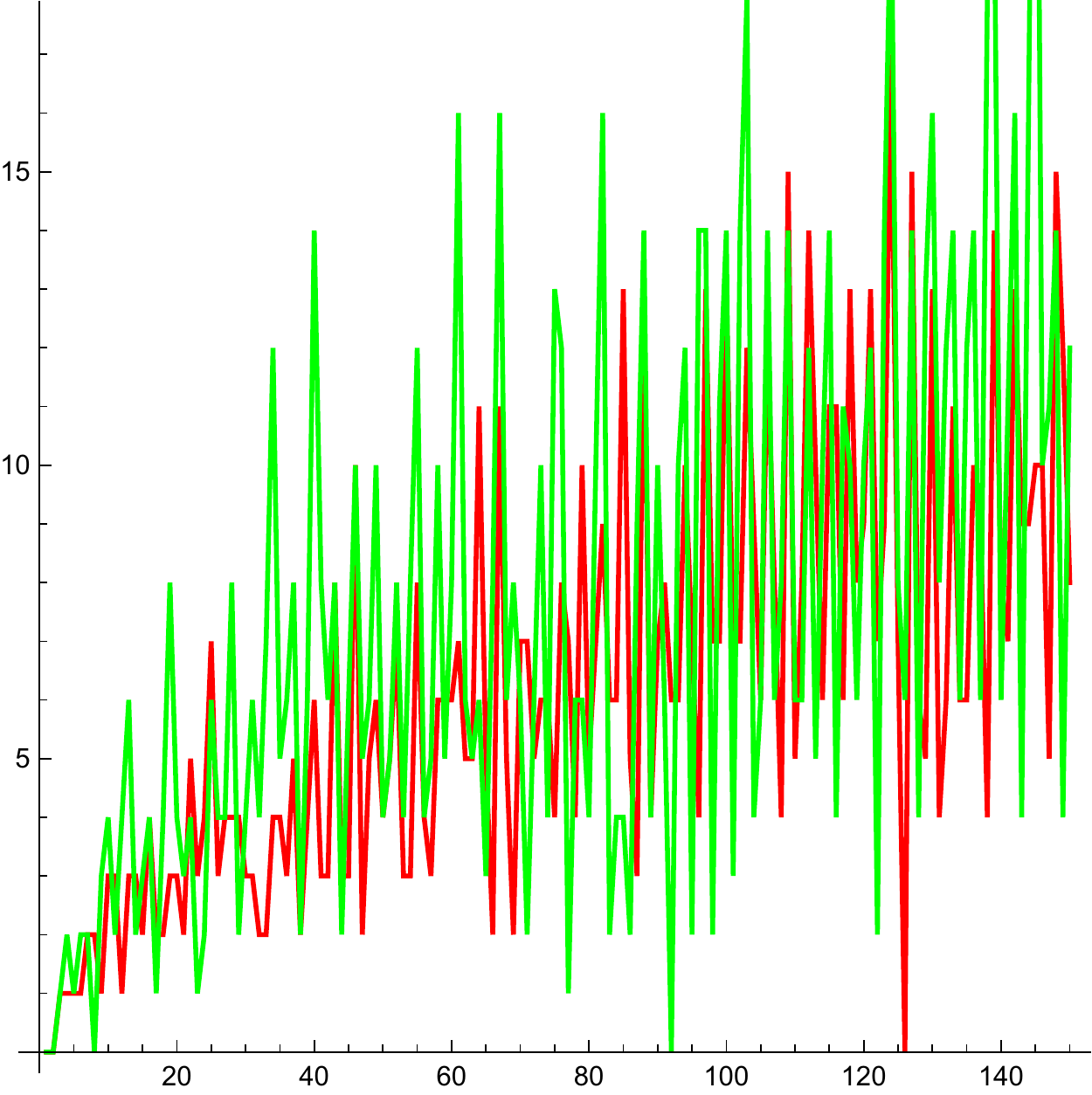}}
\scalebox{0.4}{\includegraphics{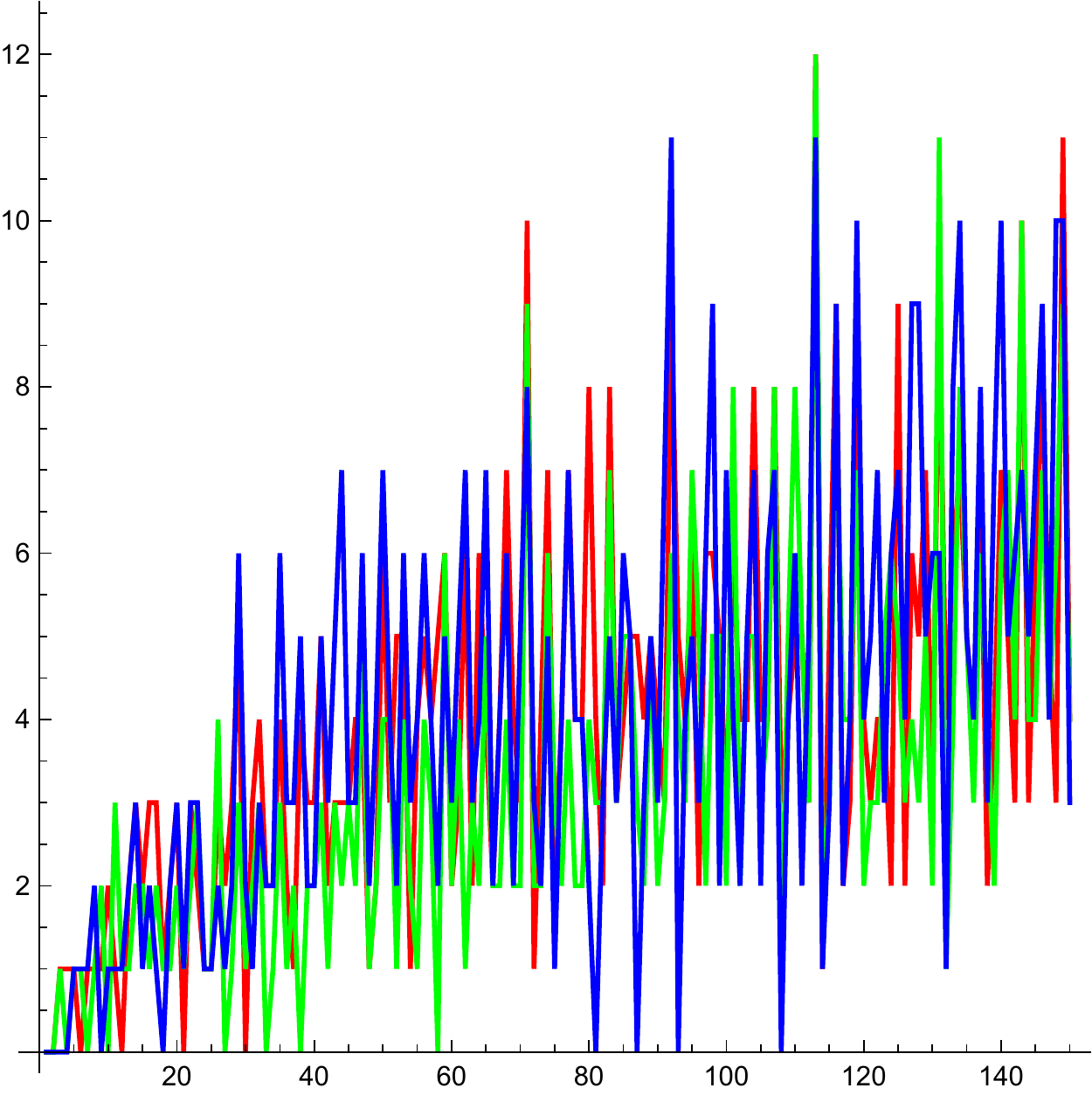}}
\caption{
The number of times an Eisenstein integer $a+ b w$ can
be written as a sum of Eisenstein primes. For $b=2$, 
this appears always possible as $p+q$ with 
form $p=x+w, q=y+w$. For $b=3$, two {\bf Eisenstein ghost twins}
which are $109+3w$ and $121+3w$ appear. Later we have still
gaps when forcing individual rows but they don't intersect. 
}
\end{figure}

\begin{figure}[!htpb]
\scalebox{0.3}{\includegraphics{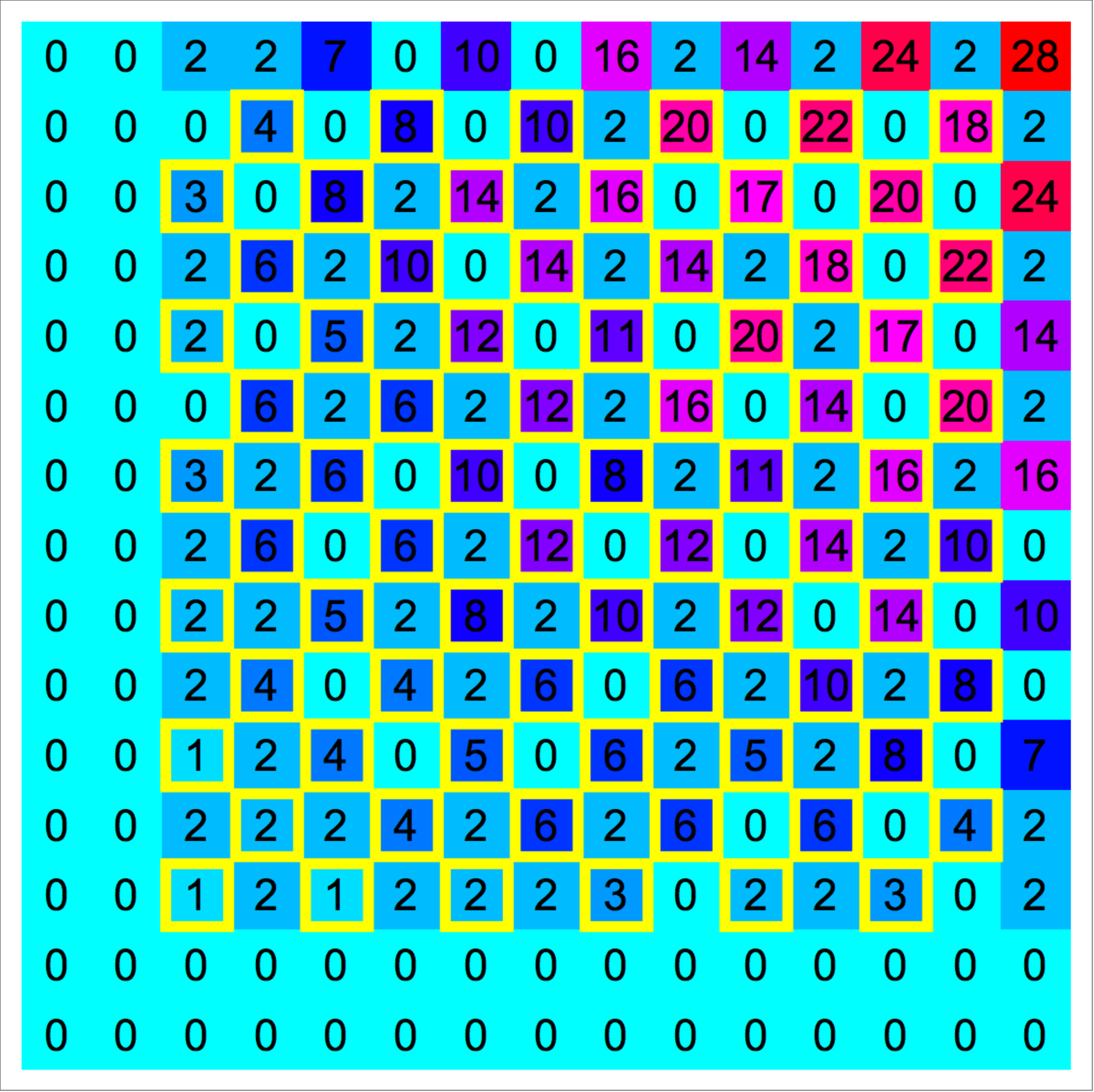}}
\caption{
The matrix entries in the picture show how many times a Gaussian integer can be written as a sum
of two Gaussian primes in $Q$. We see that the evenness condition is clearly necessary. 
\label{goldbachmatrix}
}
\end{figure}

\begin{figure}[!htpb]
\scalebox{0.3}{\includegraphics{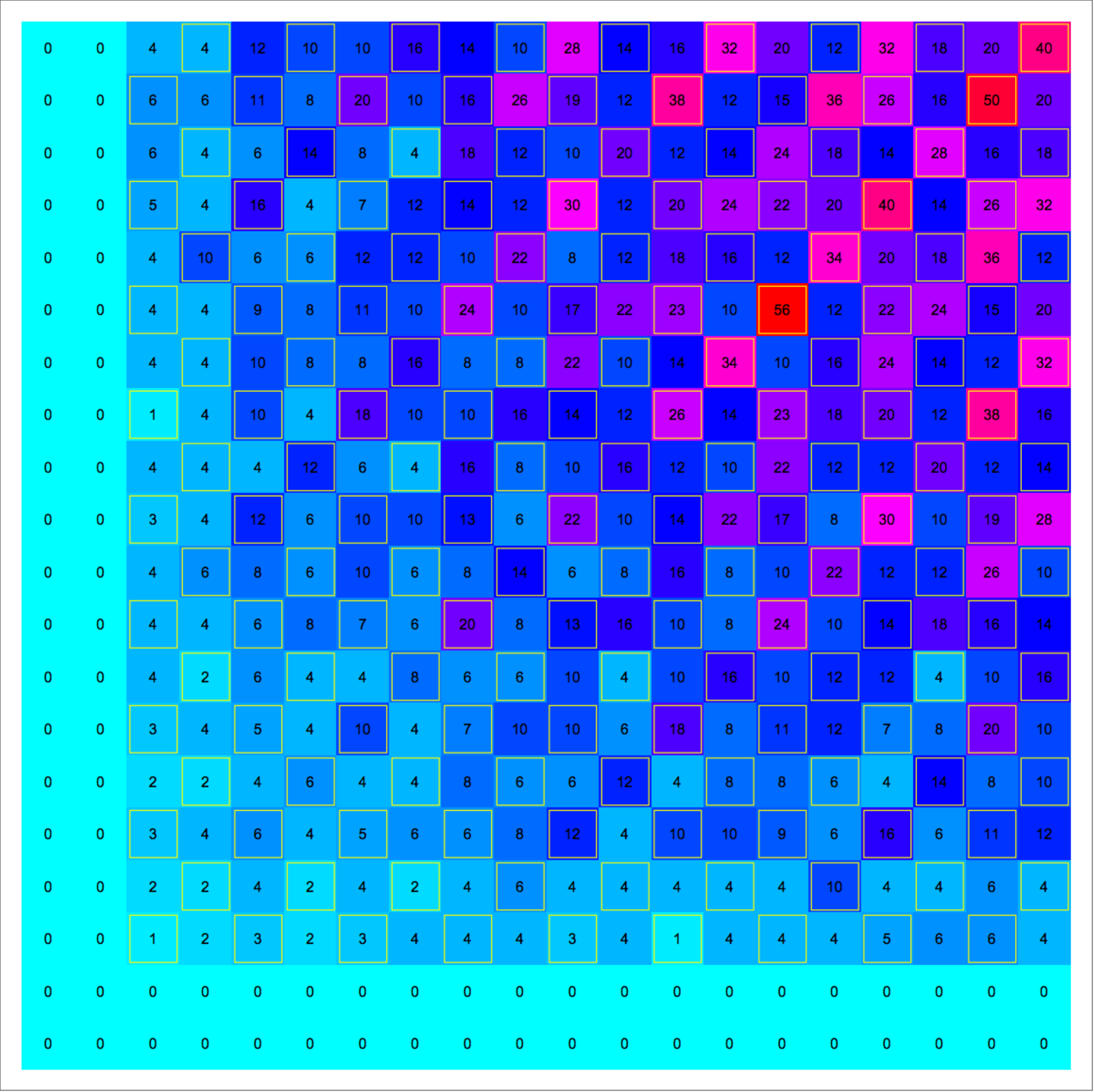}}
\caption{
Here we see how many times an Eisenstein integer $a+w b$ can be written
as a sum of two Eisenstein primes in $Q$. It looks as if we can write any $a + b w$
with $a>1,b>1$ as a sum of two Eisenstein primes in $Q$. But it only appears so at first.
}
\label{circles}
\end{figure}

\begin{figure}[!htpb]
\scalebox{0.3}{\includegraphics{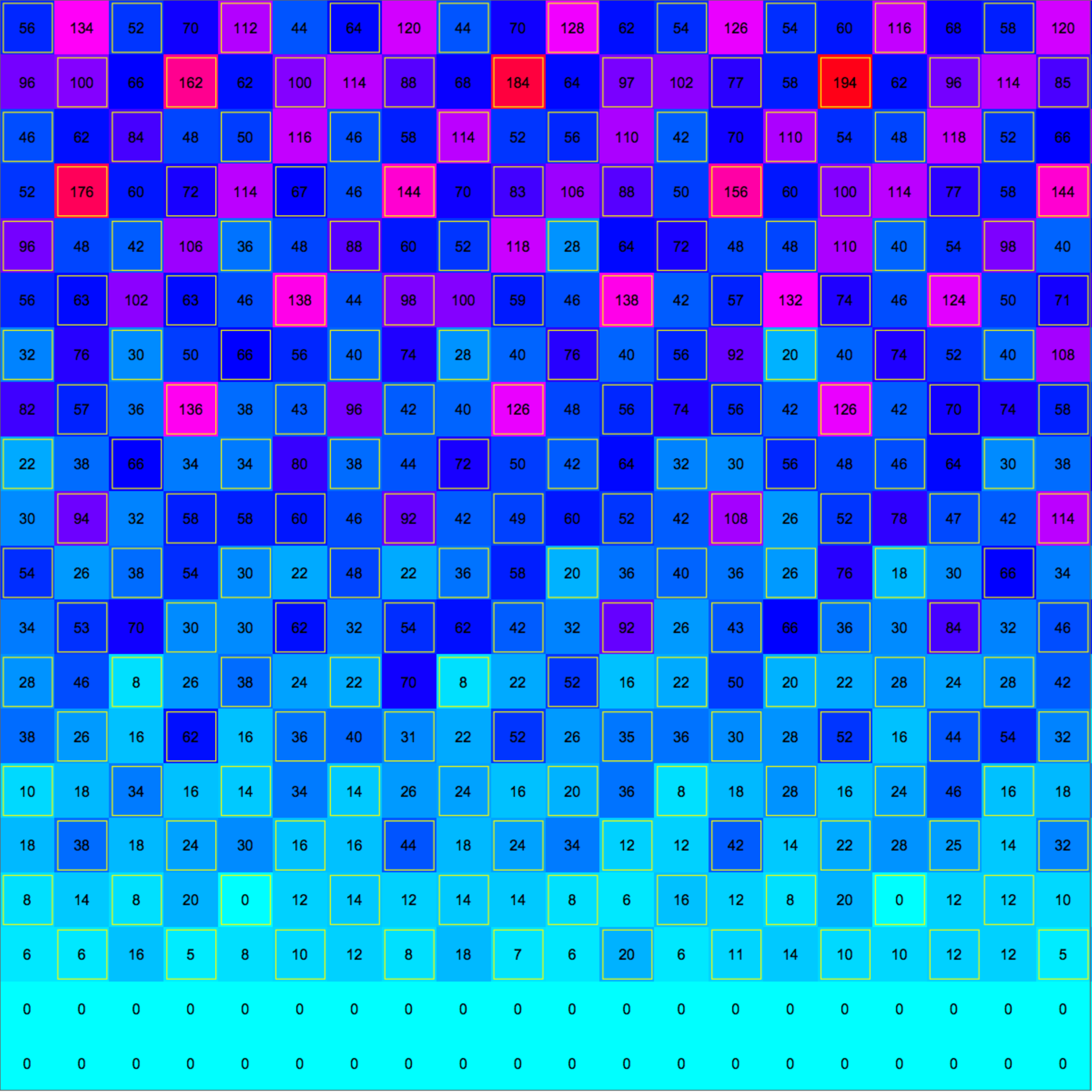}}
\caption{
There are two integers $(3+109w)$ and $(3+121w)$ (as well as their mirrors) which can not be written as the sum
of two positive Eisenstein primes $a+bw + (c+dw)$ with positive $a,b,c,d$.
These are the bad {\bf Eisenstein ghost twins}. Can you find the ghosts in 
the matrix? We believe they are the only ones. 
}
\label{circles}
\end{figure}

\bibliographystyle{plain}

\end{document}